\documentclass[twocolumn,12pt,reprint,aps,pra,longbibliography,preprintnumbers,showkeys]{revtex4-1}

\usepackage[T1]{fontenc}
\usepackage{graphicx}

\usepackage{geometry}
\usepackage{amssymb}

\usepackage{xcolor}
\usepackage{hyperref}
\newcommand{\orcid}[1]{\href{https://orcid.org/#1}{\textcolor[HTML]{A6CE39}{#1}}}

\usepackage{mathrsfs}

\begin{document}
    
\title{Rigorous computer-assisted bounds on the period doubling
  renormalisation fixed point and eigenfunctions in maps with critical
  point of degree 4}

\author{Andrew D Burbanks}
\thanks{\orcid{0000-0003-0685-6670}}
\email{andrew.burbanks@port.ac.uk}
\author{Andrew H Osbaldestin}
\thanks{\orcid{0000-0002-2482-0738}}
\author{Judi A Thurlby}
\thanks{\orcid{0000-0002-7372-9471}}
\affiliation{School of Mathematics and Physics, University of Portsmouth, UK}


%\PACS{05.10.Cc,05.45.-a}
%MSC-class: 37E20 (Primary) 37E05 (Secondary)

\date{\today}

\begin{abstract}
We gain tight rigorous bounds on the renormalisation fixed point for
period doubling in families of unimodal maps with degree $4$ critical
point.  We use a contraction mapping argument to bound essential
eigenfunctions and eigenvalues for the linearisation of the operator
and for the operator controlling the scaling of added noise.
Multi-precision arithmetic with rigorous directed rounding is used to
bound operations in a space of analytic functions yielding tight
bounds on power series coefficients and universal constants to over
320 significant figures.
\end{abstract}

\maketitle

\section{Introduction}\label{introduction}

\subsection{Background}

An explanation for the remarkable universality observed in
period-doubling cascades for families of unimodal maps of the interval
with quadratic critical point was offered by Feigenbaum
\cite{feigenbaum1978quantitative,feigenbaum1979universal,feigenbaum1979metric}
and Coullet and Tresser \cite{tresserc1978iterations} in terms of a
renormalisation operator acting on a suitable space of functions.

The explanation rests on the following conjectures: There exists a
nontrivial hyperbolic renormalisation fixed point. The spectrum of the
linearisation of the operator has a single essential expanding
eigenvalue. The associated one-dimensional unstable manifold crosses
the manifold corresponding to functions with superstable period
\(2^n\) orbits transversally for sufficiently large \(n\).

Lanford \cite{lanford1982computer} established the existence of a
nontrivial locally-unique hyperbolic fixed point of the operator by
rigorous computer-assisted means. He established that a certain
quasi-Newton operator is a contraction mapping on a carefully chosen
ball in a suitable space of functions and then bounded the spectrum of
the derivative of the operator at the fixed point in order to
establish hyperbolicity.

The efficacy of rigorous computer-assisted proofs in this area is
apparent in the body of work that followed. Eckmann et al
\cite{eckmann1982existence,eckmann1984computer} proved the existence
of a fixed point of the corresponding renormalisation operator for
period doubling in area-preserving maps, providing a detailed
framework for rigorous computation in Banach spaces of multivariate
analytic functions. Eckmann and Wittwer \cite{ew1985computer} examined
universality in period doubling for families of unimodal maps in the
limit of large even integer degree at the critical point.

These techniques have also proved effective in establishing universal
scaling results concerning the breakup of quasiperiodicity in various
scenarios. Mestel \cite{mestel1985computer} proved the existence and
hyperbolicity of a renormalisation fixed point for the breakup of
quasiperiodicity in circle maps with golden mean rotation
number. MacKay \cite{mackay1993renormalisation} examined critical
scaling in the breakup of invariant tori in area-preserving maps, and
Stirnemann \cite{stirnemann1999existence} proved the existence of the
corresponding critical fixed point for the breakup of conjagacy to
rigid rotation taking place on the boundary of Siegel discs in
iterated complex maps.

Analytical proofs of universality for critical scaling in the period
doubling of families of unimodal maps have been been harder to come
by.  Campanino et al \cite{campanino1982feigenbaum} proved existence
of the nontrivial renormalisation fixed point for period doubling in
the case of unimodal maps with degree \(2\) at the critical
point. Epstein \cite{epstein1986new} established that solutions to the
corresponding functional equation exist within the class of even
functions of general degree at the critical point providing another
proof that did not require a computer. Eckmann and Wittwer
\cite{eckmann1987complete} recast the problem in terms of an extended
renormalisation group operator, written in a form that includes the
bifurcation parameter itself, and hence established existence and
hyperbolicity of the fixed point for maps with degree \(2\) at the
critical point, together with transversal crossing of the manifold of
superstable period two functions by the corresponding unstable
manifold, thus providing a full proof of the Feigenbaum conjectures in
the case of critical exponent \(2\). The reader is referred to
Cvitanovic \cite{cvitanovic1989universality} for a thorough compendium
of results in this area.

The work of Douady and Hubbard in complexifying the operator, together
with Sullivan's program to find the fixed point
\cite{douady1985polynomial,sullivan1987quasiconformal}, enriched the
field with ideas from holomorphic dynamics, Teichmueller theory, and
hyperbolic geometry. McMullen
\cite{mcmullen1994complex,mcmullen1996fiber} developed the approach of
quasiconformal rigidity and hence established global uniqueness of the
nontrivial renormalisation fixed point. Lyubich
\cite{lyubich1999feigenbaum} and Avila and Lyubich (see, in
particular, \cite{avila2011horseshoe}) extended global uniqueness and
hyperbolicity of the fixed point to arbitrary even integer degree,
establishing the existence of a renormalisation horseshoe. Faria et al
\cite{faria2006global} have extended global hyperbolicity from
analytic to \(C^r\) mappings in the degree \(2\) case. A survey of
four decades of research in the area is provided by
\cite{lyubich2012forty}. More recently, Gorbovickis and Yampolsky
\cite{gorbovickis2018noninteger} have broadened the reach to certain
maps with non-integer critical exponent.

\subsection{Overview}

In this note, we focus on universality in period-doubling of unimodal
maps of degree \(4\) at the critical point and note that maps with
other even integer degrees are amenable to the same treatment. While
not generic, the case of degree 4 critical point may be of interest
for systems possessing certain symmetries and for the case of locally
bimodal maps in which one quadratic extremum is mapped to another.

We note that existence and hyperbolicity of the renormalisation fixed
point follows from the work of \cite{lyubich1999feigenbaum}. Our
motivation is to find tight rigorous bounds on the fixed point
function, on eigenfunctions of the derivative and of the operator
controlling the scaling of noise, and on the corresponding universal
constants.

Firstly, we use a modified operator that encodes the action of the
renormalisation operator on maps, $g$, that can be written as
$g(x)=G(x^4)$.  We adapt the methods of proof of
\cite{lanford1982computer,eckmann1982existence,eckmann1984computer,ew1985computer,mestel1985computer,mackay1993renormalisation,stirnemann1999existence},
with the addition of multi-precision arithmetic and parallel
computation. We use rigorous computer-assisted means (`function-ball
algebra') to gain tight bounds on the nontrivial fixed point of the
renormalisation operator, by showing that a quasi-Newton operator for
the fixed-point problem is a contraction map on a suitable ball in a
Banach space of analytic functions
(Sections~\ref{the-renormalisation-fixed-point},
\ref{existence-of-the-fixed-point}).

We bound the spectrum of the derivative of the operator at the fixed
point, producing crude initial bounds on the relevant eigenvalues and
establishing their multiplicities. We then take a novel approach to
bounding the eigenfunctions and eigenvalues by recasting the
eigenproblem for the derivative operator in a modified nonlinear form,
and using a contraction mapping argument. In particular, we gain tight
bounds on the eigenfunction corresponding to the essential expanding
eigenvalue delta (Section~\ref{spectral-theory}).  By adapting the
method to the relevant operator, we bound the eigenfunction and
eigenvalue that govern the universal scaling of additive uncorrelated
noise (Section~\ref{critical-scaling-of-additive-noise}).

The contraction mapping arguments used each require careful
consideration of the action of the relevant operators on high-order
terms, in order to mitigate the function-ball analogue of the
dependency problem, well-known in interval arithmetic.  We present
solutions to the corresponding dependency problems in each case.

Our computations use multi-precision arithmetic with rigorous directed
rounding modes to bound tightly the coefficients of the relevant power
series (including the polynomial parts taken to high truncation degree
alongside rigorous bounds on all high-order terms). Indeed, we are
able to obtain bounds that are tight, in the \(\ell_1\)-sense, on the
power series coefficients of the critical fixed point, on the
eigenfunctions corresponding to critical scaling in both the dynamical
space and the parameter space, and on the eigenfunction corresponding
to the critical scaling of additive noise, together with their
accompanying universal scaling constants.

Working to degree \(2560\) (reduced to \(640\) via symmetry), we are
able to bound the fixed point within a ball of analytic functions of
\(\ell_1\) radius \(10^{-331}\). We bound the eigenfunction
corresponding to the parameter-scaling eigenvalue within radius
\(10^{-325}\) and the eigenfunction controlling the scaling of
additive noise within radius \(10^{-323}\). We note that the
individual power series coefficients of these functions are therefore
constrained within intervals having those same radii. This yields
bounds on universal scaling constants in both the dynamical and the
parameter space, and on the eigenvalue for scaling of additive noise:
we are able to prove \(331\), \(325\), and \(323\), digits of these
correct, respectively, confirming and extending significantly the
accuracy of previous numerical estimates.

\section{The renormalisation fixed point}\label{the-renormalisation-fixed-point}

\subsection{The renormalisation operator}\label{the-renormalisation-operator}

We consider the operator \(R\) defined by
\begin{equation}
Rg(x):= a^{-1}g(g(ax)),\label{eqn:R}
\end{equation}
where \(a=a_g:= g(1)\) is chosen to preserve the normalisation
\(g(0)=1\). (We note that other choices for \(a\), also preserving
this normalisation, may be taken and that, as is well-known, the
particular variant will later affect the spectrum of \(DR(g)\) only up
to coordinate-change eigenvalues.)

We seek a nontrivial fixed point of \(R\), with a critical point of
degree $4$ at the origin, in a Banach space \(A:=\mathscr{A}(\Omega)\)
of functions analytic on an open disc
\(\Omega=D(c,r):=\{z\in\mathbb{C}:\ |z-c|<r\}\) and continuous on its
closure, \(\overline{\Omega}\), with (finite) \(\ell_1\)-norm.

We work with a modified operator that encodes the action of $R$ on
functions possessing the symmetry $g(x)=G(x^4)$. Specifically, we let
\(X = Q(x) := x^4\) and write
\[
g(x) = G(Q(x)) = G(X). 
\]
We then seek a fixed point of the corresponding operator \(T\)
defined by
\begin{equation}
TG(X) := a^{-1}G(Q(G(Q(a)X))),\label{eqn:T}
\end{equation}
where \(a:= G(1)\).

\subsection{The disc algebra}\label{the-disc-algebra}

We write \(G\in\mathscr{A}(\Omega)\) as
\[
G = G^{u}\circ\psi,
\]
where
\(\psi:\overline{\Omega}\to\overline{\mathbb{D}}:=\overline{D(0,1)}\)
is the affine map from the domain \(\overline{\Omega}\) to the unit
disc given by
\[
\psi:x\mapsto\frac{x-c}{r}.
\]
We then take \(G^{u}\in \mathscr{A}(\mathbb{D})\), the disc algebra:
the set of functions analytic on the open unit disc \(\mathbb{D}\) and
continuous on its closure, \(\overline{\mathbb{D}}\), with (finite)
\(\ell^1\)-norm. Equipped with the usual addition and scalar
multiplication, viz. \((f+g)(x)=f(x)+g(x)\) and \((b f)(x)=b f(x)\),
and with the \(\ell^1\)-norm, \(\mathscr{A}(\mathbb{D})\) (and, hence,
\(\mathscr{A}(\Omega\))) is a Banach space (moreover, when equipped
with the product \((f\cdot g)(x)=f(x)\cdot g(x)\), it is a commutative
unital Banach algebra) isometrically isomorphic to the sequence space
\(\ell^1\); functions \(f\in \mathscr{A}(\Omega)\) may be written as
power series expansions
\[
f(x) = \sum_{k=0}^{\infty} a_k\left(\frac{x-c}{r}\right)^k,
\]
convergent on \(\Omega\).

\subsection{Nonrigorous calculation}\label{nonrigorous-calculation}

Firstly, we compute approximate fixed points of the renormalisation
operator, \(T\), by working {in the space of power series truncated to
  some fixed degree} \(N\) expanded on the disc \(\Omega\). To this
end, we write \(\ell_1\) as the direct sum,
\[
\ell_1\cong\mathbb{R}^{N+1}\oplus\ell_1,
\]
and let \(PA\) and \(HA=(I-P)A\) denote the canonical projections onto
the polynomial part and high-order part of the space, respectively.
Thus we may write \(f\in A\) as
\[
f = f_P + f_H,
\]
with \(f_H\in HA\) and \(f_P\in PA\) where
\[
f_P(x) = \sum_{k=0}^{N} a_k\left(\frac{x-c}{r}\right)^k.
\]
As a starting point, we consider the one-parameter family of maps
given by
\[
f_\mu(x) = 1-\mu x^4,
\]
and choose a parameter value \(\mu\) close to the accumulation
\(\mu_\infty\) of the first period-doubling cascade for the
family. (The intention is to find a function that lies close to the
stable manifold of the critical renormalisation fixed point.) We
establish, {by making use of multi-precision arithmetic to locate}
superstable periodic orbits of periods \(2^k\) for \(1\le k\le 32\),
that \(\mu_\infty\simeq 1.594901356228820564497828\). Writing
\(f_\mu=G\circ Q\) and then applying the renormalisation operator
iteratively until we no longer observe an improvement in the residue
\(\|T^{n+1}(G)-T^{n}(G)\|\) (when working with our chosen truncation
degree and precision) then provides an initial approximate fixed
point.

\subsection{Newton operator}\label{newton-operator}

\begin{figure}[ht]\begin{center}
\includegraphics[width=0.425\textwidth]{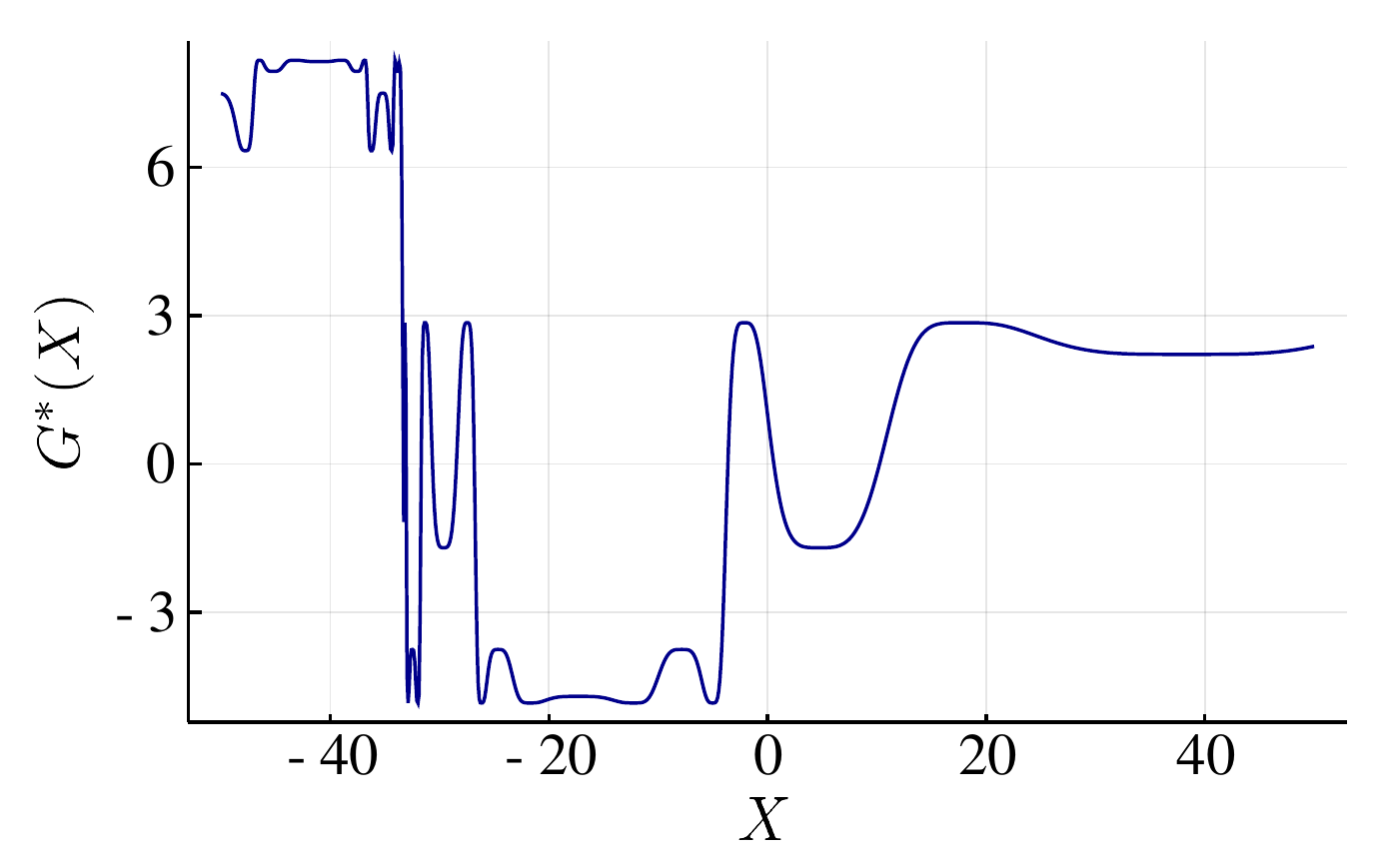}\\
\includegraphics[width=0.425\textwidth]{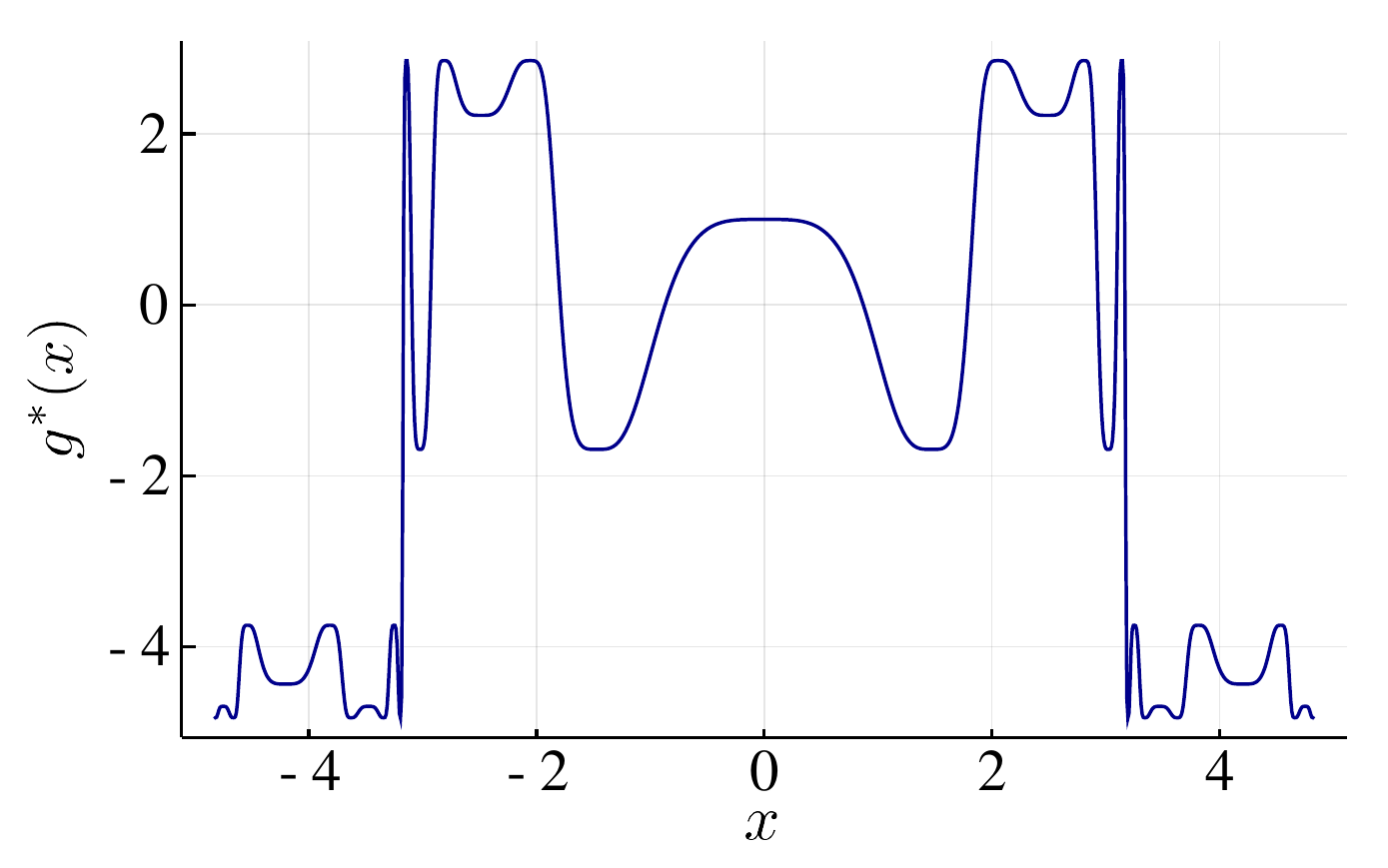}\\
\includegraphics[width=0.425\textwidth]{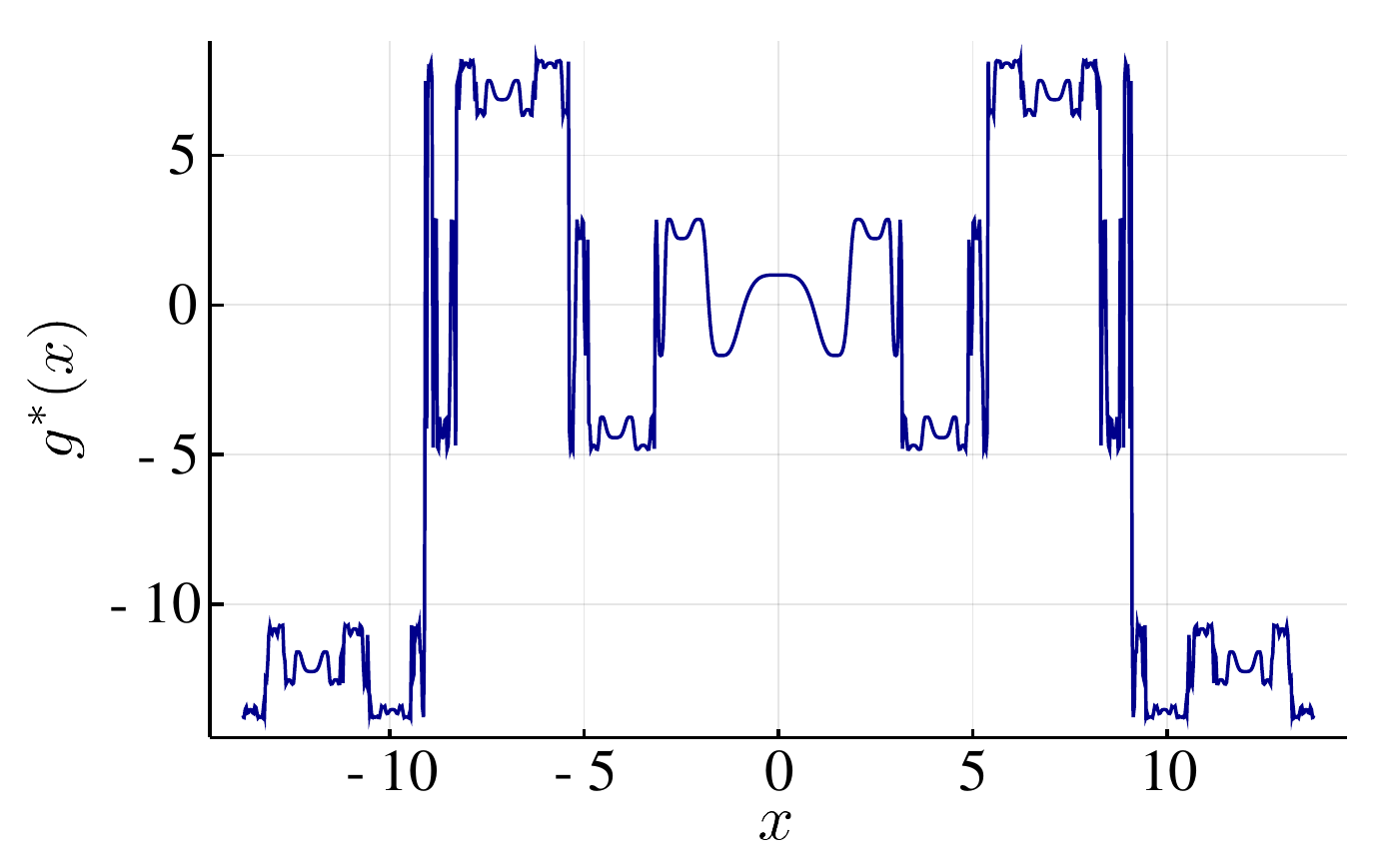}\\
\includegraphics[width=0.425\textwidth]{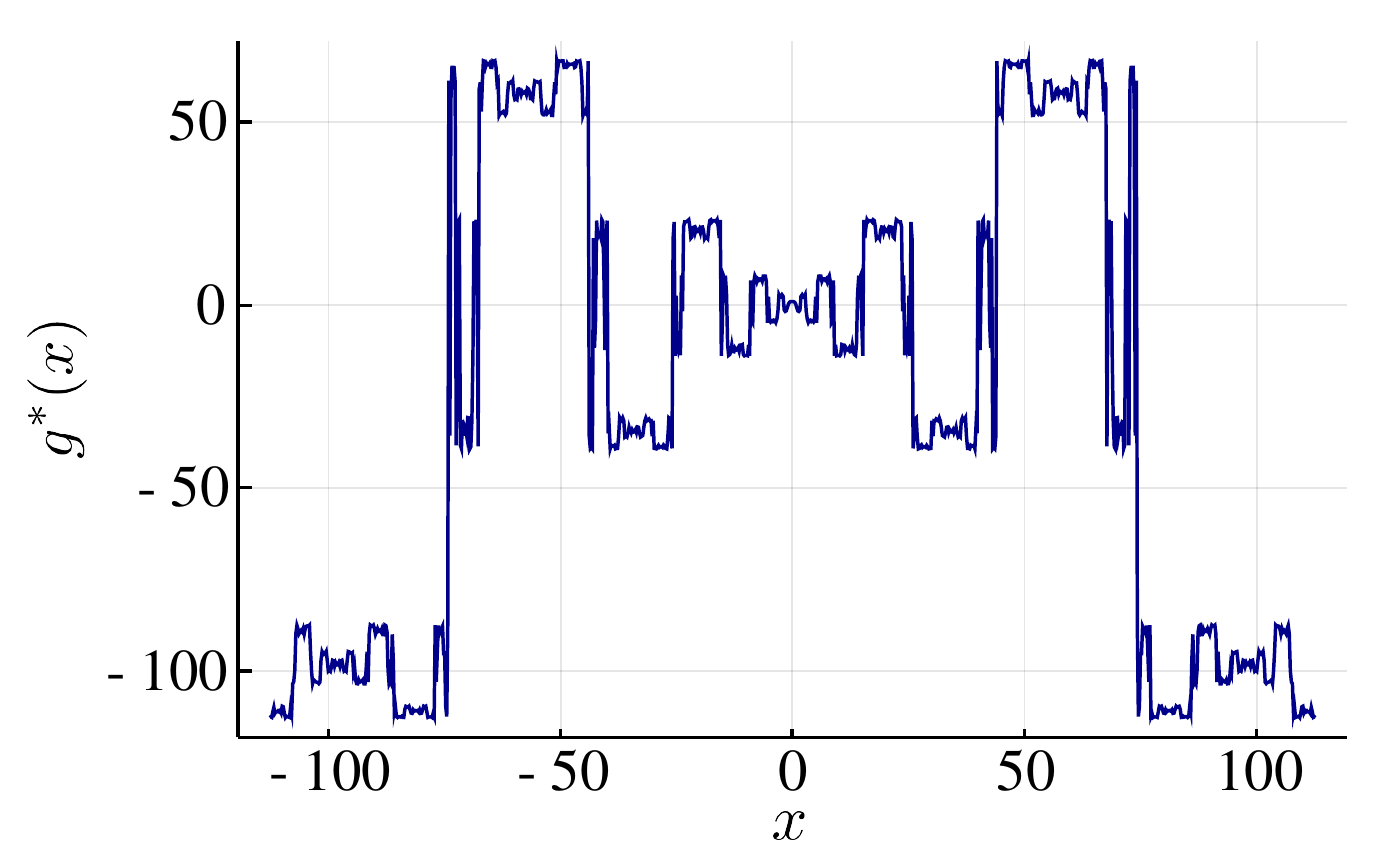}
\caption{Nonrigorous approximations of the function $G^{*}$ (topmost) and $g^{*}(z)$ for $z\in[-|\alpha|^{k},|\alpha|^{k}]$ for $k=3, 5, 9$, where $\alpha=1/a$ and $a=G^{*}(1)$.\label{fig:gn}}
\end{center}\end{figure}

We note that fixed points of \(T\) are zeros of the operator \(F =
T-I\), and perform Newton iterations, in the space of power series
truncated to degree \(N\), to approximate such a zero. The one-step
Newton operator is given by
\begin{eqnarray}
\phi:G
&\mapsto& G - [DF(G)]^{-1}F(G)\nonumber\\
&=& G - [DT(G) - I]^{-1}(T(G) - G),\label{eqn:newton}
\end{eqnarray}
in which \(DT(G)\in\mathscr{B}(A, A)\) denotes the tangent map of
\(T\) at \(G\), given formally by the Fr\'echet derivative
\begin{widetext}
   
\begin{eqnarray}
DT(G):\delta G
&\mapsto&
-a^{-2}\delta a G(Q(G(Q(a)X)))\nonumber\\
&&{}+ a^{-1}\biggl(
    \delta G(Q(G(Q(a)X)))\nonumber\\
&&\quad{}+ G'(Q(G(Q(a)X)))\cdot Q'(G(Q(a)X))\cdot \bigl[
        \delta G(Q(a)X)\nonumber\\
&&\qquad{}+ G'(Q(a)X)\cdot Q'(a)\delta a\cdot X
    \bigr]
\biggr),\label{eqn:frechet}
\end{eqnarray}
\end{widetext}
where \(\delta a = \delta G(1)\). After the Newton iterations converge
to our chosen precision, we denote the resulting approximate fixed
point by \(G^0\). (See Fig. \ref{fig:gn}.)

Our goal is then to appeal to the contraction mapping theorem to prove
that the operator \(T\) has a locally-unique fixed point in a ball
\(B^1\) of functions centered on \(G^0\) in the space
\(\mathscr{A}(\Omega)\). The operator \(T\) is not itself contractive
at the fixed point (indeed, we later bound the spectrum of the
derivative there and obtain the eigenfunctions corresponding to the
expanding eigenvalues). However, we can find a quasi-Newton operator
\(\Phi\) that has the same fixed points as \(T\) and establish instead
that \(\Phi\) is a contraction mapping on \(B^1\).

\section{Bounds on the fixed point}\label{existence-of-the-fixed-point}

\subsection{Rigorous computations in the function space}\label{rigorous-computations-in-the-function-space}

We bound operations in the function space \(\mathscr{A}(\mathbb{D})\)
(and hence \(\mathscr{A}(\Omega)\)) by maintaining careful control
over the coefficients of truncated power series along with all
high-order terms. In order to maintain rigour, we work with interval
arithmetic using high-precision computer-representable bounds with
directed rounding modes, conforming to the relevant industry standards
(section~\ref{sec:computational-issues}).  To this end, we define a
ball of functions, centered on a polynomial \(f_P\in
P\mathscr{A}(\mathbb{D})\), with high-order bound \(v_H\ge 0\) and
general ``error'' bound \(v_E\ge 0\), as follows
\begin{eqnarray*}
B(f_P; v_H, v_E)
&:=& \bigl\{
f\in \mathscr{A}(\mathbb{D}):\\
&&\ f = f_P+f_H+f_E,\\
&&\ f_H\in H\mathscr{A}(\mathbb{D}),\|f_H\|\le v_H,\\
&&\ f_E\in \mathscr{A}(\mathbb{D}),\|f_E\|\le v_E
\bigr\}.
\end{eqnarray*}

Following \cite{eckmann1984computer,ew1985computer}, we extend the
definition slightly, to the case where the function \(f_P\) is not
known exactly, but rather has coefficients \(a_k\) confined within
intervals.  Let \(v_P = ([b_0,c_0],\ldots,[b_N,c_N])\in J^{N+1}\) be a
vector of intervals (here, \(J\) denotes
\(\{[a,b]:\ a,b\in\mathbb{R},a\le b\}\)).  Given the bounds
\(v=(v_P,v_H,v_E)\), we define the standard function ball
\(B(v)\subset\mathscr{A}(\mathbb{D})\) by
\begin{eqnarray*}
B(v_P, v_H, v_E)
&:=&
\bigl\{
f\in\mathscr{A}(\mathbb{D}):\\
&&\ f = f_P+f_H+f_E,\\
&&\ f_P\in P\mathscr{A}(\mathbb{D}),\\
&&\ f_P(x) = \sum_{k=0}^{N}a_kx^k,a_k\in[b_k,c_k],\\
&&\ f_H\in H\mathscr{A}(\mathbb{D}),\|f_H\|\le v_H,\\
&&\ f_E\in \mathscr{A}(\mathbb{D}),\|f_E\|\le v_E
\bigr\}.
\end{eqnarray*}
The resulting set of functions is convex and closed. The definition
extends in a natural way to function balls \(B_\Omega(v_P,v_H,v_E)\),
for a general disc \(\Omega\), by writing \(f=f^{u}\circ\psi\) where
\(f^{u}\in B(v_P,v_H,v_E)\).

We bound operations on the function space \(\mathscr{A}(\Omega)\) by
first choosing computer-representable numbers for the quantities in
\(v:=(v_P,v_H,v_E)\). For each binary operation \(\oplus\), we then
implement a version, \(\oplus_b\), acting on bounds \(v,w\) such that
for all \(f\in{B}(v)\) and \(g\in{B}(w)\),
\[
f \oplus g\in{B}(v)\oplus{B}(w)\subseteq{B}(v\oplus_b w).
\]
The operation \(v\oplus_bw\), on bounds, is constructed carefully in
order to guarantee that the above inclusion holds even when
implemented using finite-precision arithmetic. In this way, all vector
space operations, together with the product, \(f\cdot g\), composition
of functions, \(f\circ g\), differentiation followed by composition,
\(f'\circ g\), and the norm \(\|f\|\), may be bounded. For an
exhaustive exposition, in the case of maps of two variables, see
\cite{eckmann1984computer}.

\subsection{Quasi-Newton operator}\label{quasi-newton-operator}

The Newton operator for the fixed-point problem was shown in equation
\ref{eqn:newton}. However, in order to establish contractivity, we
would need to work with its derivative, which would involve taking the
second Fr\'echet derivative of \(T\). This proves to be inconvenient
in practice. Instead, we note that if \(\Lambda\) is any invertible
linear operator, then the fixed points of the quasi-Newton method
given by
\begin{equation}
\Phi:G \mapsto G - \Lambda(T(G)-G),\label{eqn:quasinewton}
\end{equation}
are exactly the fixed points of \(T\). We choose
\[
\Lambda \simeq [DT(G)-I]^{-1},
\]
and establish that our chosen $\Lambda$ is indeed invertible.
Specifically, we approximate the Fr\'echet derivative \(DT(G)\) by a
fixed linear operator \(\Delta\simeq DT(G^0)\) with action zero on
high-order terms. For the polynomial terms, we evaluate the expression
for the Fr\'echet derivative at Schauder basis elements given by the
sequence of monomials
\[
e_j(x) = \left(\frac{x-c}{r}\right)^j,
\]
for \(j=0,\ldots,N\) and bound the resulting matrix elements by
trivial intervals to give a real interval matrix denoted
\(\Delta_{PP}\).  We compute an interval matrix \(\Lambda_{PP}\)
guaranteed to bound the inverse \((\Delta_{PP}-I)^{-1}\). Thus the
corresponding linear operator \(\Lambda\) has action \(\Lambda_{PP}\)
on the polynomial part of the space, and action \(-I\) on the
high-order part.

\subsection{Bound 1: distance moved by the approximate fixed point}\label{bound-1-distance-moved-by-the-approximate-fixed-point}\label{sec:epsilon}

In order to use the contraction mapping theorem, we need to prove that
a certain ball in \(\mathscr{A}(\Omega)\) is mapped into itself
contractively by \(\Phi\). We achieve this by establishing two bounds:
a bound on how far the approximate fixed point \(G^0\) moves under the
operator \(\Phi\), and a bound on the derivative \(D\Phi\) that we
will use in order to show that \(\Phi\) is contractive and that
\(\Phi\) maps the ball to itself.

To this end, we define a ball of functions \(B^0=B_\Omega(G^0;0,0)\)
of radius zero; the singleton \(\{G^0\}\). By applying \(\Phi\) to
\(B^0\), in the sense of using the corresponding function ball
operations to find a new function ball guaranteed to contain the
result, we gain a rigorous bound on how far \(G^0\) moves under
\(\Phi\):
\begin{equation}
\|\Phi(G^0)-G^0\| < \varepsilon.\label{eqn:epsilon}
\end{equation}
We now choose a radius \(\rho > \varepsilon\) and form the function
ball \(B^1=B_\Omega(G^0;0,\rho)\), on which we need to prove that
\(\Phi\) is a contraction mapping.

\subsection{Domain extension}\label{sec:domain-extension}

\begin{figure}[ht]\begin{center}
\includegraphics[width=0.425\textwidth]{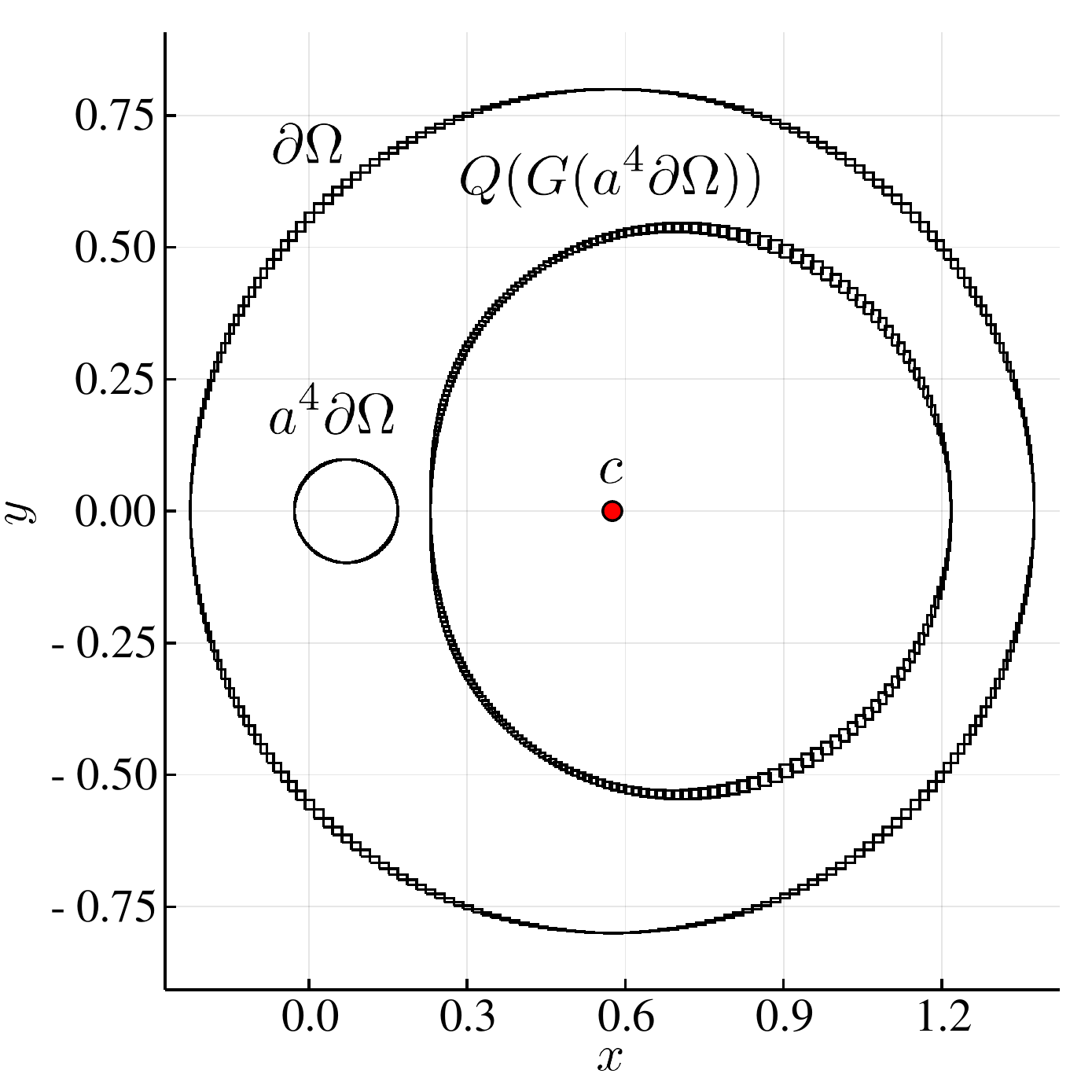}
\caption{Verification of the domain extension conditions computed using the function ball $B^1$, using a rigorous covering of the boundary $\partial\Omega$ by $256$ rectangles.\label{fig:de}}
\end{center}\end{figure}

The first step in what follows is to show that \(T\) is well-defined
and differentiable, with compact derivative, on \(B^1\). We do this by
establishing the `domain extension' or `analyticity improving'
property \cite{mackay1993renormalisation,mestel1985computer}: for all
\(G\in B^1\) we demand that
\begin{eqnarray}
Q(a)\overline{\Omega}
&\subset& \Omega\label{eqn:de1},\\
Q(G(Q(a)\overline{\Omega}))
&\subset& \Omega\label{eqn:de2}.
\end{eqnarray}
In the above, the overline denotes topological closure.  Recall that
we take \(a:= G(1)\). Thus the universal quantifier is not vacuous for
equation \ref{eqn:de1}. Systematic experimentation is used to find a
suitable domain \(\Omega=D(c,r)\). For the case of quartic critical
point, we may choose \(\Omega=D(0.5754, 0.8)\). The domain may be
improved further by choosing \(c\) so as to minimise the absolute
value of the constant term on \(G^u\) where
\(G^0=G^u\circ\psi\). Doing so reduces the dominant contribution to
the error bounds involved in composition. Fig.  \ref{fig:de}
illustrates domain extension for a rigorous covering of the boundary
\(\partial\Omega\). Rectangle arithmetic, in which intervals bound the
real and imaginary parts of rectangles covering $\partial\Omega$, is
used to confirm the result.

We note that the space \(\mathscr{A}(\Omega)\) is infinite-dimensional
and has the bounded approximation property
\cite{kato2013perturbation,yosida2013functional}. It follows that the
spectrum of a compact operator consists of \(0\) together with only
isolated eigenvalues of finite multiplicity. The spectrum of
finite-rank approximations converges to the spectrum of the operator
itself; if \(L\) is compact and \(\|L'-L\|\to 0\), then the spectrum
of \(L'\) (and, indeed, the corresponding eigenfunctions) converges to
that of \(L\) apart from at \(0\)
\cite{krasnoselskij1972approximate}. (In the case of complex domains,
one can prove that domain extension yields compactness by appealing to
the Cauchy estimates on suitable discs to provide uniform continuity,
and hence establish normality. Montel's theorem then implies the
result \cite{mackay1993renormalisation}.) Compactness will prove
crucial in bounding the spectrum of the linearisation \(DT(G)\) at the
fixed point in section \ref{sec:spectrum}.

\subsection{Bound 2: uniform contractivity}\label{bound-2-uniform-contractivity}\label{sec:kappa}

Our final goal is to find a uniform bound on the contractivity of
\(\Phi\) on \(B^1\). We do this by bounding
\begin{equation}
\|D\Phi(G)\| \le \kappa < 1,\quad\mbox{for all}\ G\in B^1,\label{eqn:kappa}
\end{equation}
for a suitable norm, and then appealing to the mean value
theorem (that this yields uniform contractivity may be seen by
considering the line segment joining any two points in the convex set
\(B^1\) and noting that a bound on the norm of \(D\Phi(G)\) valid for
all \(G\in B^1\) provides an upper bound on all of the corresponding
pairwise contractivities).

The Fr\'echet derivative of the quasi-Newton operator \(\Phi\) (from
equation \ref{eqn:quasinewton}) is given by
\begin{equation}
D\Phi(G): \delta G
\mapsto \delta G-\Lambda[DT(G)\delta G-\delta G].\label{eqn:dphi}
\end{equation}

We bound \(D\Phi(G)\) {via} the maximum column sum norm. That is, we
bound the norms \(\|D\Phi(G)e_k\|\) for all basis elements \(e_k\) and
then take the supremum, noting that
\[
\|D\Phi(G)\|
:= \sup_{\|\delta G\|=1}\|D\Phi(G)\delta G\|
\le \sup_{k}\|D\Phi(G)e_k\|,
\]
where the norm on the left is the standard operator norm.  To do this,
we bound the action of the Fr\'echet derivative of \(\Phi\) at \(B^1\)
on function balls containing the \(e_k\). Firstly, we let \(E_k :=
B(e_k;0,0)\) for \(k=0,\ldots,N\), i.e., we consider singletons
containing each of the polynomial basis elements. The problem of
capturing the (infinitely-many) norms that remain is reduced to a
finite computation by taking the single ball \(E_H := B(0;1,0)\),
i.e., the convex hull of all high-order basis elements, and bounding
\(\|D\Phi(B^1)E_H\|\), i.e., \(\|D\Phi(G)\delta G\|\) for all \(G\in
B^1\) and \(\delta G\in E_H\). This yields
\[
\kappa\ge\sup\left\{\|D\Phi(G)E\|\right\},
\]
with the supremum taken over all \(G\in B^1\) and all
\(E\in\{E_0,\ldots,E_N,E_H\}\), from which, for \(\kappa<1\), the mean
value theorem delivers the uniform bound on contractivity
\[
\|\Phi(f)-\Phi(g)\| \le \kappa\|f-g\|\quad\mbox{for all}\ f,g\in B^1.
\]

\begin{figure}[ht]\begin{center}
\includegraphics[width=0.45\textwidth]{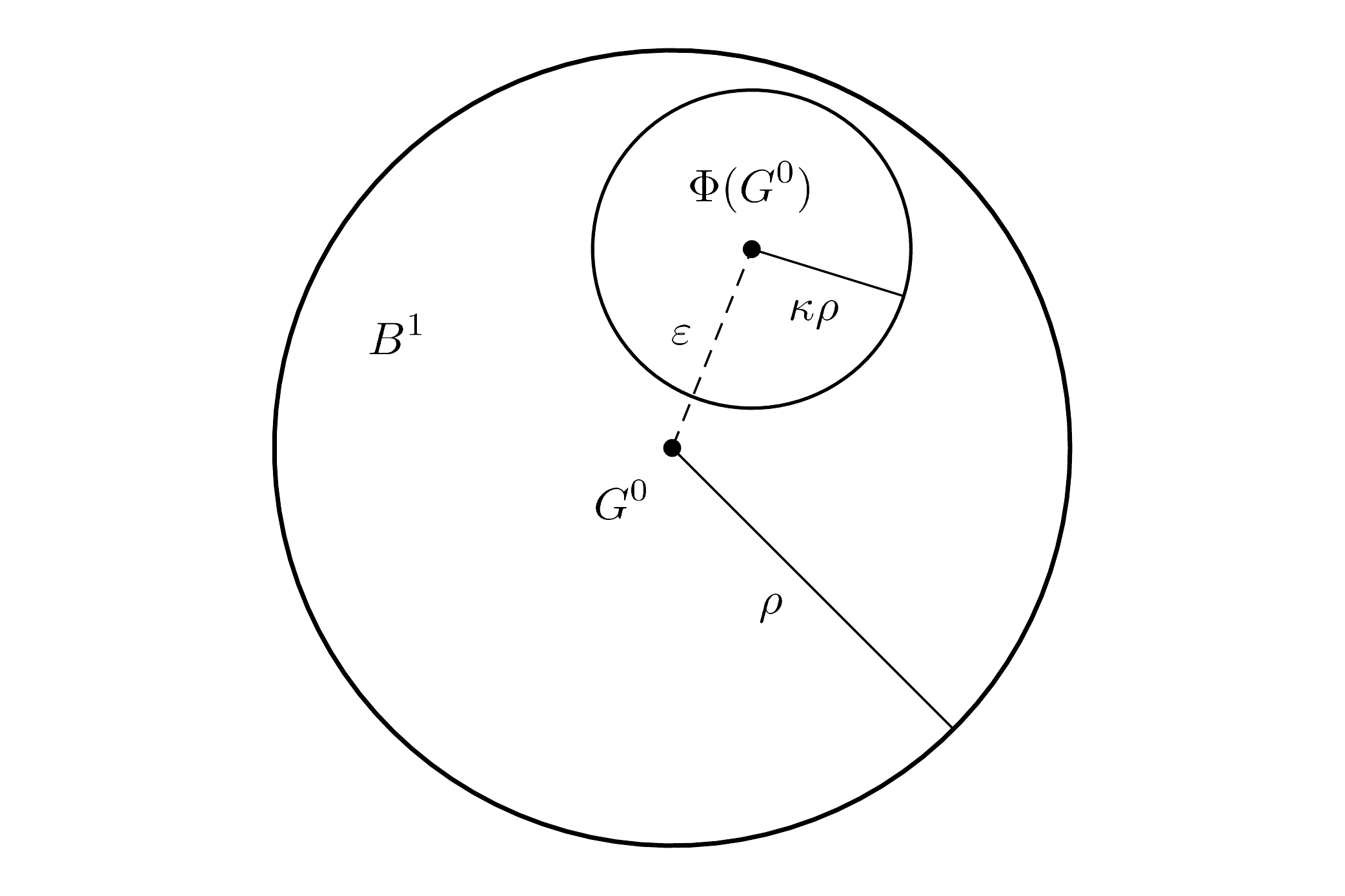}
\caption{Schematic of the contraction mapping.\label{fig:cmt}}
\end{center}\end{figure}

\subsection{Dependency problems}\label{dependency-problems}\label{sec:depend}

It is crucial, for the case where \(\delta G\) is a high-order
perturbation, i.e., \(\delta G\in H\mathscr{A}(\Omega)\), to mitigate
the function ball analogue of the dependency problem, well-known in
interval arithmetic \cite{moore1966interval,kaucher2014self}. In the
expression for \(D\Phi\) (equation \ref{eqn:dphi}), the action of
\(\Lambda\) on high-order terms is \(-I\), thus the action of
\(D\Phi\) on a high-order perturbation \(\delta b_H\) is given by:
\begin{eqnarray}
D\Phi(B^1)\delta b_H
&=& \delta b_H - \Lambda\left[DT(B^1)\delta b_H - \delta b_H\right]\label{eqn:depend1}\\
&=& \delta b_H - \Lambda\left[DT(B^1)\delta b_H\right] - \delta b_H\label{eqn:depend2}\\
&=& -\Lambda\left[DT(B^1)\delta b_H\right].\label{eqn:depend3}
\end{eqnarray}
Computing the norm \(\|D\Phi(B^1)E_H\|\) naively by performing
function ball operations based on expression \ref{eqn:depend1} would
result in an upper bound on contractivity larger than \(2\), even in
the case where \(D\Phi(B^1)\) is indeed contractive, due to the
implicit presence of uncancelled terms \(\delta b_H-\delta b_H\) in
\ref{eqn:depend2}. The operands in an expression of the form
\(\|f-g\|\) {where the subtraction operation is implemented in
  function-ball arithmetic} are treated as independent (high-order)
functions, here, subject only to the bounds \(\|f\|,\|g\|\le
1\). Expression \ref{eqn:depend3} must therefore be used instead.

\subsection{Existence and local uniqueness}\label{existence-and-local-uniqueness}

Finally, using the bounds obtained in equations \ref{eqn:epsilon} and
\ref{eqn:kappa}, we verify the inequality
\[
\varepsilon < \rho(1-\kappa),
\]
to ensure that \(\Phi(B^1)\subset B^1\), which establishes that
\(\Phi\) is a contraction mapping on \(B^1\). Fig. \ref{fig:cmt}
illustrates the situation schematically. Hence, \(\Phi\) (and,
therefore, \(T\)) has a locally unique fixed point, \(G^{*}\in B^1\).

Using our chosen disc \(\Omega\), we are able to complete the proof by
choosing truncation degree \(N=40\), thus \(g\) has degree
\(160\). Working with precision equivalent to \(40\) digits in the
significand, we obtain \(\varepsilon=1.59\times 10^{-21}\), and
choosing \(\rho=10^{-20}\) gives \(\kappa=6.88\times 10^{-3}\).

\begin{figure}[ht]\begin{center}
\includegraphics[width=0.425\textwidth]{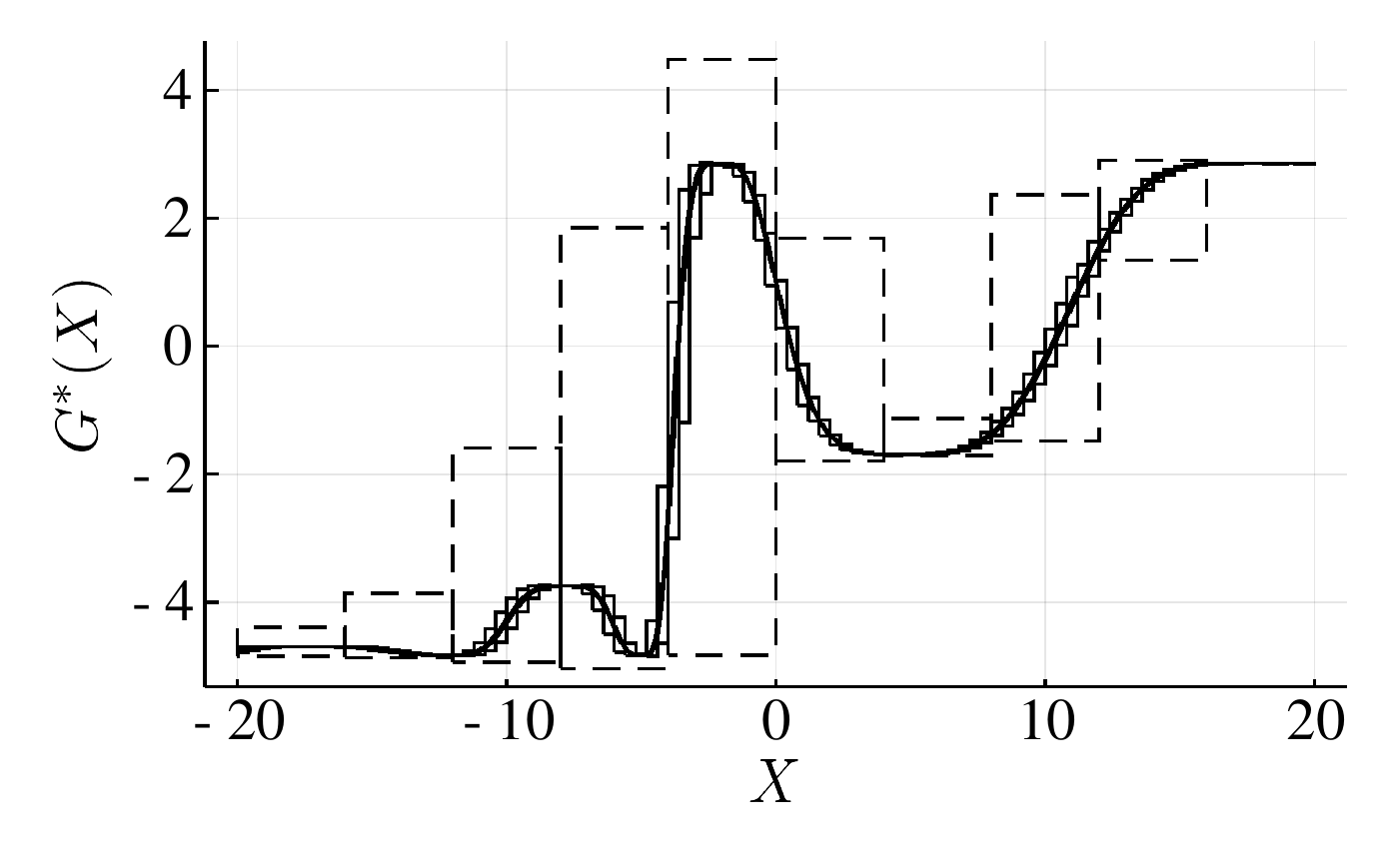}\\
\includegraphics[width=0.425\textwidth]{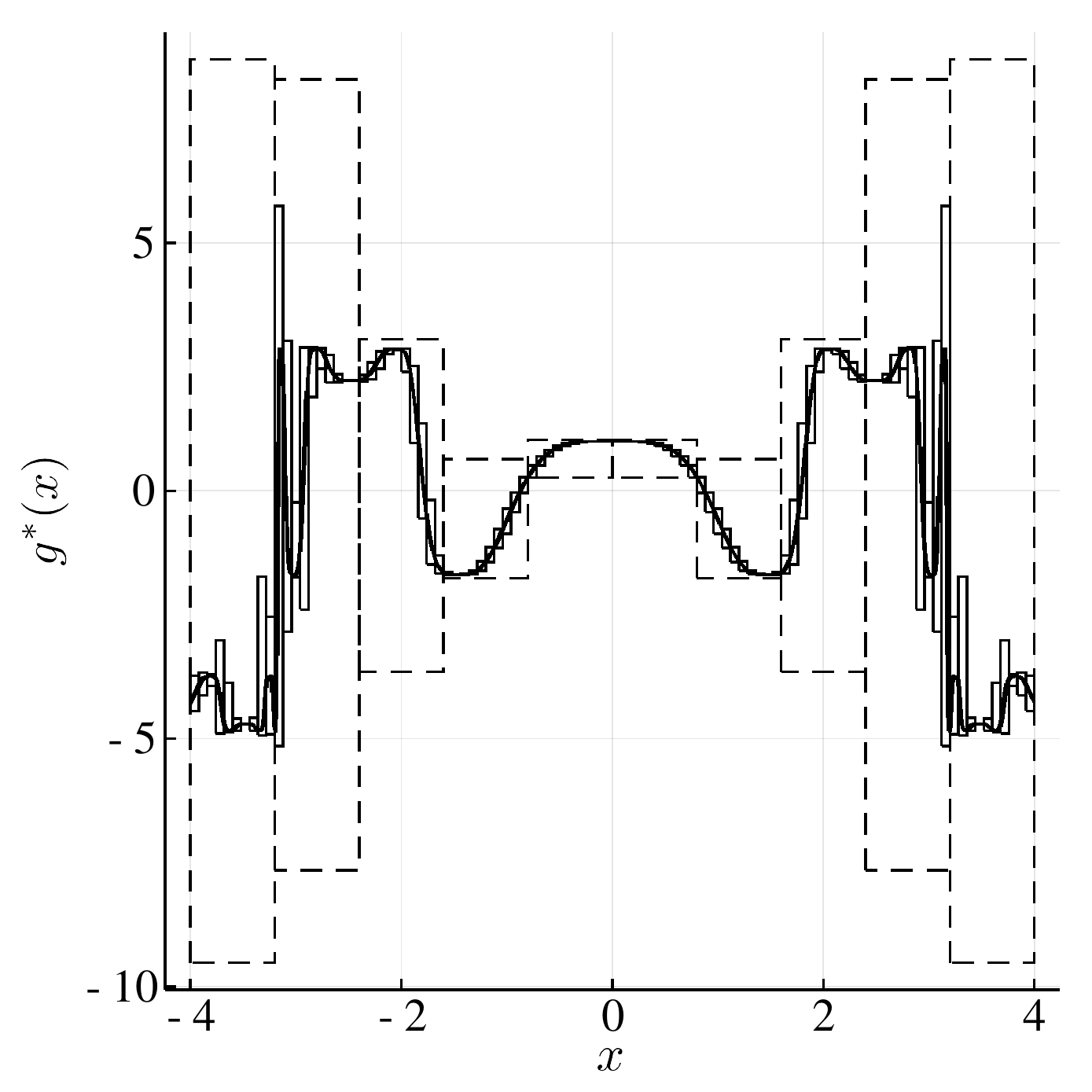}
\caption{Rigorous coverings of the functions $G^{*}$ (top) and $g^{*}$
  (bottom) using 10 (dashed lines), 50, and 500 rectangles computed
  using the function ball $B^1\ni G^{*}$ together with the fixed-point
  equation.\label{fig:coverg}}
\end{center}\end{figure}

\begin{table*}[ht]
{\small\begin{tabular}{|r|r||l|l|l||l|l|l||l|l|l|}
\hline

\multicolumn{2}{|c||}{} &
\multicolumn{3}{c||}{Fixed point ($G^{*}$)} &
\multicolumn{3}{c||}{Delta eigenfunction ($V^{*}$)} &
\multicolumn{3}{c|}{Noise eigenfunction ($W^{*}$)} \\

\hline

$N$ & $\#bits$ &
$\varepsilon$ & $\rho$ & $\kappa$ &
$\hat\varepsilon$ & $\hat\rho$ & $\hat\kappa$ &
$\tilde\varepsilon$ & $\tilde\rho$ & $\tilde\kappa$ \\

\hline\hline

$40$ & $132$ &
$1.59\cdot 10^{-21}$ & $10^{-20}$ & $6.88\cdot 10^{-3}$ &
$3.17\cdot 10^{-16}$ & $10^{-15}$ & $1.17\cdot 10^{-3}$ &
$2.35\cdot 10^{-16}$ & $10^{-15}$ & $7.85\cdot 10^{-3}$ \\

$80$ & $265$ &
$3.75\cdot 10^{-42}$ & $10^{-41}$ & $1.01\cdot 10^{-6}$ &
$4.88\cdot 10^{-37}$ & $10^{-36}$ & $1.39\cdot 10^{-7}$ &
$7.33\cdot 10^{-37}$ & $10^{-36}$ & $1.90\cdot 10^{-7}$ \\

$160$ & $531$ &
$7.84\cdot 10^{-84}$ & $10^{-83}$ & $1.36\cdot 10^{-12}$ &
$8.37\cdot 10^{-79}$ & $10^{-78}$ & $1.87\cdot 10^{-13}$ &
$8.57\cdot 10^{-78}$ & $10^{-77}$ & $4.32\cdot 10^{-14}$ \\

$320$ & $1063$ &
$2.89\cdot 10^{-166}$ & $10^{-165}$ & $3.01\cdot 10^{-24}$ &
$1.52\cdot 10^{-160}$ & $10^{-159}$ & $4.12\cdot 10^{-25}$ &
$6.24\cdot 10^{-160}$ & $10^{-159}$ & $9.56\cdot 10^{-26}$ \\

$480$ & $1594$ &
$4.14\cdot 10^{-249}$ & $10^{-248}$ & $7.28\cdot 10^{-36}$ &
$2.21\cdot 10^{-243}$ & $10^{-242}$ & $9.99\cdot 10^{-37}$ &
$5.31\cdot 10^{-242}$ & $10^{-241}$ & $2.32\cdot 10^{-37}$ \\

$640$ & $2126$ &
$5.01\cdot 10^{-332}$ & $10^{-331}$ & $1.85\cdot 10^{-47}$ &
$2.90\cdot 10^{-326}$ & $10^{-325}$ & $2.53\cdot 10^{-48}$ &
$1.36\cdot 10^{-324}$ & $10^{-323}$ & $5.87\cdot 10^{-49}$ \\

\hline
\end{tabular}}
\caption{Parameters for {contraction mappings and resulting bounds on}
  the renormalisation fixed point, $G^{*}$, the eigenfunction,
  $V^{*}$, corresponding to $\delta$, and the eigenfunction, $W^{*}$,
  corresponding to the scaling of additive noise.  In all cases, the
  number of digits $P$ in the significand, for the decimal
  floating-point versions of the proofs, was chosen to be equal to the
  truncation degree $N$.  (The table also indicates the corresponding
  number of bits chosen in the significand for the independent binary
  floating point versions of the proofs.  Experimentation reveals that
  we may reduce $P$ at least as far as $\left\lfloor
  2N/3\right\rfloor$, for the computations shown, and still gain
  rigorous bounds of the same orders of magnitude.)\label{tab:bounds}}
\end{table*}

\begin{table*}[ht]
{\small\begin{tabular}{rllllllll}
$a={}$
-0.&5916099166 &3443815013 &9624354381 &6289537902
   &2298919075 &5829639056 &2608082701 &6110024444\\
   &6553096873 &1159671843 &1035214180 &0643269743
   &8637238931 &2068288207 &7993159616 &2409259411\\
   &5430529642 &7613470988 &2939926870 &4915779588
   &8740837617 &0145437404 &8090852176 &8119211417\\
   &0711171042 &5330824210 &0970358064 &2260084834
   &3287080164 &7846778564 &3980486155 &4138928900\\
   &8050440114 &\ldots     &           &          
   & & & & \\
$\alpha=1/a={}$
-1.&6903029714 &0524485334 &3780150324 &1613482282
   &7805970956 &1966682423 &2634497392 &1908881055\\
   &1432766085 &7861529191 &5193152630 &8212594164
   &1050775616 &3090857294 &0573192526 &2783102042\\
   &4401895602 &5177655047 &9352262368 &7664454132
   &1907107192 &6768349355 &4697194567 &2766866785\\
   &1484514531 &8901391119 &4135568528 &2120804754
   &6969604755 &8987391859 &3295066623 &5922528661\\
   &8546743362 &\ldots     &           &
   & & & & \\
$\delta={}$
+7.&2846862170 &7334336430 &8930567995 &5530694780
   &4661979979 &0659072121 &2901883462 &1435067620\\
   &0657264503 &1360371147 &0784357866 &9255573693
   &3221121594 &9170167056 &0272610414 &2834709598\\
   &2287873290 &2387885867 &2064166568 &1895073101
   &1658106317 &3127916581 &6323366267 &7746542527\\
   &7844194832 &0362437902 &4983698686 &8146702404
   &9663158059 &7051641021 &9527093166 &3172744588\\
   &9929\ldots &           &           &          
   & & & & \\
$\gamma={}$
+8.&2439108542 &5258681839 &8462365029 &2376160673
   &1776662405 &8409262192 &5682565366 &3924142562\\
   &6899642047 &2075784242 &2300873689 &8322349635
   &1071732825 &3743947119 &1666888923 &2401827811\\
   &4543435570 &5947708003 &7798523831 &6683467659
   &8572907048 &7598764245 &8476648182 &5677074055\\
   &9568984297 &6849327088 &1184491967 &8812146275
   &7670908015 &1177052580 &3233041606 &2789993350\\
   &21\ldots   &           &           &          
   & & & &
\end{tabular}
}
\caption{Digits proven correct of $a=G^{*}(1)$ ($331$ digits),
  $\alpha=1/a$ ($331$ digits), $\delta=\varphi(V^{*})$ ($325$ digits),
  and $\gamma=\varphi(W^{*})$ ($323$ digits) obtained from the proof
  with truncation degree $N=640$ for $G^{*},V^{*},W^{*}$
  (corresponding to degree $4N=2560$ for
  $g^{*},v^{*},w^{*}$).\label{tab:digits}}
\end{table*}

\subsection{Tight bounds on the fixed point}\label{tight-bounds-on-the-fixed-point}

If the goal were to provide an alternative (computer-assisted) proof
of existence of the fixed point, a relatively low truncation degree
for \(G\) and a relatively low precision is adequate (indeed, one
could even have used standard \(64\)-bit double precision numbers,
with careful control over directed rounding modes). The resulting
function ball radius \(\rho\) gives an \(\ell^1\)-bound on the
accuracy of the intervals bounding the coefficients of \(G\).

However, we improve these bounds significantly by both increasing the
truncation degree and by using rigorous multi-precision
arithmetic. Table \ref{tab:bounds} shows parameters and bounds proven
valid for establising the existence of \(G^{*}\) and hence
\(g^{*}\). Table \ref{tab:digits} lists the digits of the relevant
universal constants, including \(a:= g^{*}(1)\), that we have been
able to prove correct as a result. In particular, we prove $331$
significant digits of $a$ and $\alpha=1/a$ correct (for comparable
numerical estimates, see \cite{briggs1991precise,briggs1998analytic}).

Figure \ref{fig:coverg} demonstrates a rigorous covering of the
fixed-point functions \(G^{*}\) (resp. \(g^{*}\)). These were computed
by using the function ball \(B^1\) (resp. \(B^1\circ Q\)) with
truncation degree \(40\) on the domain \(\Omega\) (resp. on the
preimage \(Q^{-1}(\Omega)\)) together with recurrences derived from
the corresponding fixed-point equations in order to bound the
functions on larger subsets of \(\mathbb{R}\setminus\Omega\) (resp. on
its preimage under \(Q\)).  That the relevant analytic extensions
exist is ensured by the domain extension property verified in
section~\ref{sec:domain-extension}.

\section{Spectral theory}\label{spectral-theory}

\subsection{The spectrum}

Our goal here is to gain tight rigorous bounds on the spectrum of the
derivative $DT(G)$ at the fixed point, and on the corresponding
eigenfunctions. The space \(A\) is infinite-dimensional and has the
approximation property. Thus, compactness of a bounded linear operator
\(L\in\mathscr{B}(A, A)\) implies that the spectrum of \(L\) consists
of the origin together with a countable set of isolated eigenvalues of
finite multiplicity (which accumulate at \(0\))
\cite{kato2013perturbation}.

We note that the spectrum of \(DT(G)\) and that of \(DR(g)\) are
related in the following manner. Consider \(G\in A\) and \(\delta G\in
A\) and let \(g=G\circ Q\) and \(\delta g=\delta G\circ Q\), then we
have
\[
(DT(G)\delta G)\circ Q = DR(g)\delta g.
\]
Then \(\lambda\in\sigma(DT(G))\) with \(DT(G)V=\lambda V\) implies
that \(\lambda\in\sigma(DR(g))\) with \(DR(g)v=\lambda v\) where
\(v=V\circ Q\).

The spectrum of \(DT(G^{*})\) has \(2\) eigenvalues (each of
multiplicity $1$) in the complement of the closed unit disc,
\[
\alpha^4,\delta,
\]
whereas the spectrum of \(DR(g^{*})\) has \(5\) eigenvalues in the
complement of the closed unit disc,
\[
\alpha^4,\delta,\alpha^3,\alpha^2,\alpha,
\]
(the latter three correspond to perturbations ruled-out for \(DT(G)\)
on symmetry grounds) with the others in the open unit disc. Note that
\(\alpha^4\) is a coordinate-change eigenvalue and that
\(\alpha^3,\alpha^2,\alpha^1\) correspond to perturbations that
destroy the symmetry of the quartic critical point. The eigenvalue
\(\alpha^2\) plays a role in tricritical vector scaling for locally
bimodal maps in which one quadratic extremum is mapped to another,
corresponding to an additional solution
\(q_2(x)=g^{*}(\sqrt{x})^2=G^{*}(x^2)^2\) of the functional equation
\(R(g)=g\) with universal scaling constant \(g^{*}(1)^2=\alpha^2\)
\cite{chang1981tricritical,fraser1984vector}. We note also that the
choice of a particular normalisation fixing \(g(0)=1\) affects the
spectrum only up to coordinate-change eigenvalues.

\subsection{Establishing hyperbolicity}\label{establishing-hyperbolicity}

We are interested, here, in bounding eigenvalues (crudely), and
establishing their multiplicities, for the purpose of matching with
eigenfunction-eigenvalue pairs on which we will gain much tighter
bounds, below, by using a novel method with a modified eigenproblem.

Apart from non-essential eiegnvalues, the only part of the spectrum of
\(DT(G^{*})\) outside the unit disc is the eigenvalue \(\delta\)
associated with critical scaling in the parameter space for the
period-doubling cascade. All other eigenvalues are contained in the
interior of the unit disc. Subject to a projection removing
coordinate-change directions and their corresponding eigenvalues, this
helps to establish the picture conjectured by Feigenbaum, in which
\(G^{*}\) has essentially a one-dimensional unstable manifold with
eigenvalue \(\delta\) and co-dimension one stable manifold. Transverse
intersection of the unstable manifold with manifolds of superstable
periodic functions has been established elsewhere. The conclusion is
that families of maps with critical point of degree \(4\) that exhibit
a period doubling cascade, and so (generically) cross the stable
manifold transversally, display an asymptotically self-similar
bifurcation diagram with accumulation rate of period doublings given
by \(\delta\).

\begin{figure}[ht]\begin{center}
\includegraphics[width=0.425\textwidth]{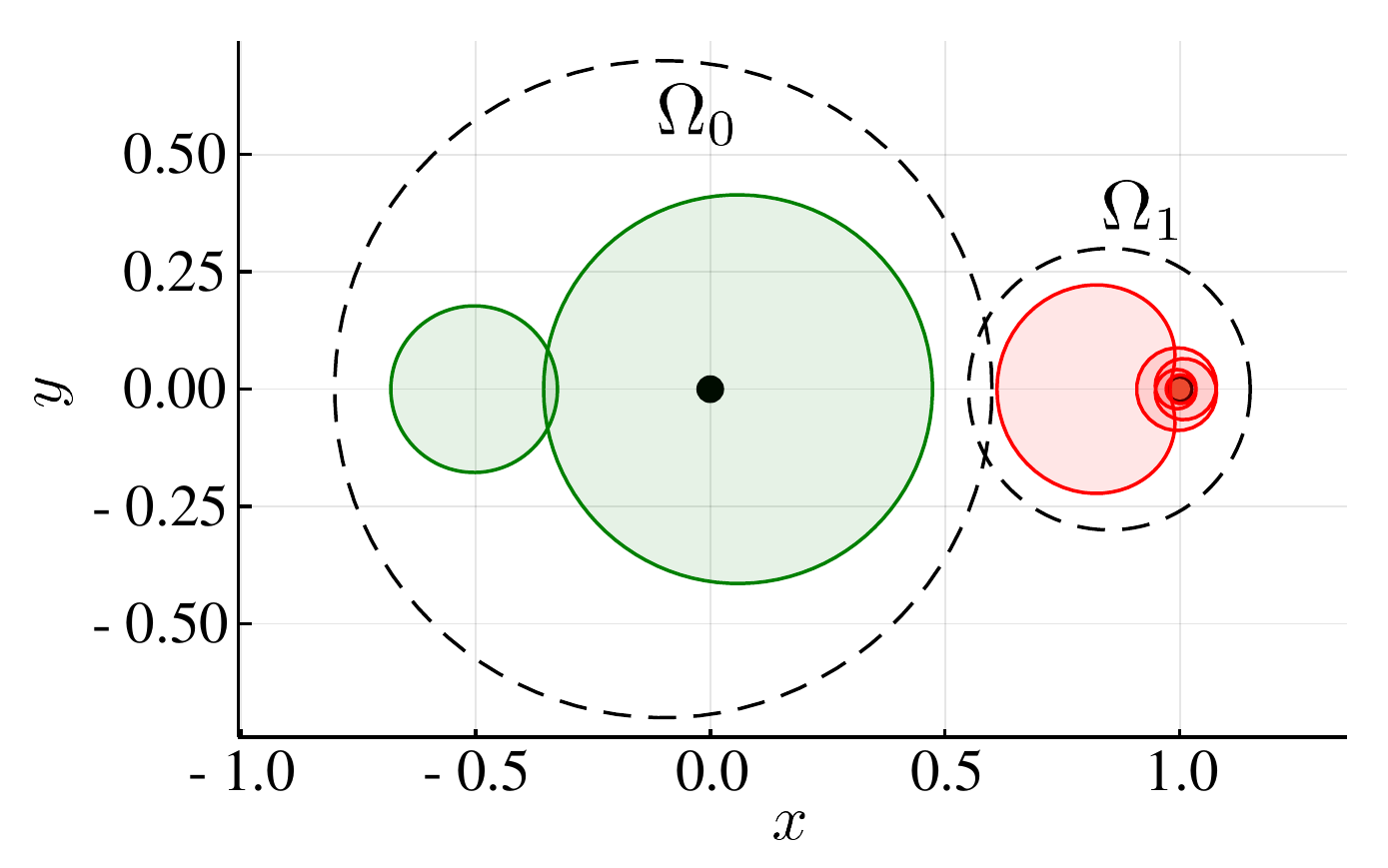}
\caption{Domain extension for $R(g)$ working in the space of pairs;
  $g=g_0\oplus g_1$ defined on domain $\Omega=\Omega_0\cup\Omega_1$
  (dashed lines), showing that $a\overline{\Omega}\subset\Omega_0$ (on
  the left; green in colour copy) and
  $g(a\overline{\Omega})\subset\Omega_1$ (on the right; red in colour
  copy).\label{fig:depair}}
\end{center}\end{figure}

We bound the spectrum for \(DT(G^{*})\) and also for \(DR(g^{*})\)
directly. We first outline the differences for \(DR(g)\) before
presenting the method common to both. The Fr\'echet derivative of
\(R\) is given formally by
\begin{eqnarray}
DR(f)\delta f
&=& -a^{-2}\delta a f(f(ax))\nonumber\\
&&{}+ a^{-1}\delta f(f(ax))\label{eqn:drg2}\\
&&{}+ a^{-1}f'(f(ax))\delta f(ax)\label{eqn:drg3}\\
&&{}+ a^{-1}f'(f(ax))f'(ax)\delta a x,\nonumber
\end{eqnarray}
where \(\delta a=\delta f(1)\).

In order to define a suitable space of functions in which to work with
\(R\), we require a domain \(\Omega\) for \(g\), with
\(0,1\in\Omega\), that satisfies the correponding domain extension
conditions
\begin{eqnarray}
a\overline{\Omega}
&\subset& \Omega,\\
g(a\overline{\Omega})
&\subset& \Omega.
\end{eqnarray}
In the quartic case to hand, no single disc that works could be found.
However, it is possible to find a union of two discs that is suitable.
Thus, when working with \(R\) and \(DR(g)\), we represent \(g\) by a
pair of power series; let \(g=g_0\oplus
g_1\in\mathscr{A}(\Omega_0)\times\mathscr{A}(\Omega_1)\) with domain
\(\Omega=\Omega_0\cup\Omega_1\) where \(\Omega_0=\mathbb{D}(c_0,r_0)\)
and \(\Omega_1=\mathbb{D}(c_1,r_1)\) with
\(0\in\Omega_0,1\in\Omega_1\) and
\(\Omega_0\cap\Omega_1\neq\emptyset\). We obtain a Banach space by
choosing a norm
\[
\|g\|=\|g_0\|+\|g_1\|,
\]
corresponding to an \(\ell^1\)-norm on
\(\mathscr{A}(\Omega_0)\times\mathscr{A}(\Omega_1)\cong
\ell^1\oplus\ell^1\).  The corresponding domain maps are
\(\psi_0,\psi_1\), where \(\psi_k:x\mapsto (x-c_k)/r_k\). The power
series that we work with are therefore those for \(g_0^u\oplus
g_1^u\in\mathscr{A}(\mathbb{D}(0,1))^2\), where
\(g_k=g_k^u\circ\psi_k\).

Choosing, for example, \(\Omega=\mathbb{D}(-0.1,
0.7)\cup\mathbb{D}(0.85, 0.3)\) and noting that, in the operator, we
have \(a:= g(1)=g_1(1)\), we are able to prove that
\begin{eqnarray}
a\overline{\Omega}
&\subset& \Omega_0,\\
g(a\overline{\Omega})
&\subset& \Omega_1,
\end{eqnarray}
which yields domain extension (Fig. \ref{fig:depair}); thus \(R\) is
well-defined on the resulting space, differentiable, and the
derivative is compact.

We may complete the contraction mapping proof for \(R\) directly by
using a ball around an approximate fixed point in the space of pairs
of maps and, by choosing a suitable basis for the space, we may then
bound the spectrum of \(DR(g)\) at the fixed point directly, allowing
perturbations that destroy the symmetry \(g=G\circ Q\) (albeit at the
cost of working in the space of pairs of maps).

\subsection{Bounding the spectrum}\label{bounding-the-spectrum}\label{sec:spectrum}

We establish firstly that the spectrum has the form described above
and gain initial bounds on the eigenvalues. For brevity, we
demonstrate this for \(DT(G^{*})\) (and apply a similar procedure
directly to \(DR(g^{*})\)). To do this, we make an invertible change
of coordinates that puts \(DT(G)\) into a form \(C^{-1}DT(G)C\) close
to diagonal, for all \(G\in B^1\). We then bound the resulting
operator by a so-called contracted matrix \(M\). This is an
\((m+1)\times(m+1)\) matrix of rectangles,
\([a,b]+i[c,d]\subset\mathbb{C}\), with \(m\le N\) with the property
that if \(\lambda=[e,f]+i[g,h]\subset\mathbb{C}\) is a rectangle
containing an eigenvalue of \(C^{-1}DT(G)C\), then taking the
determinant \(\det (M-\lambda I)\) using rectangle arithmetic (a
natural complex analogue of interval arithmetic) yields a rectangle
containing zero. Thus, if the determinant is bounded away from zero,
then we conclude that the rectangle \(\lambda\) does not contain an
eigenvalue.

\begin{figure}[ht]\begin{center}
\includegraphics[width=0.425\textwidth]{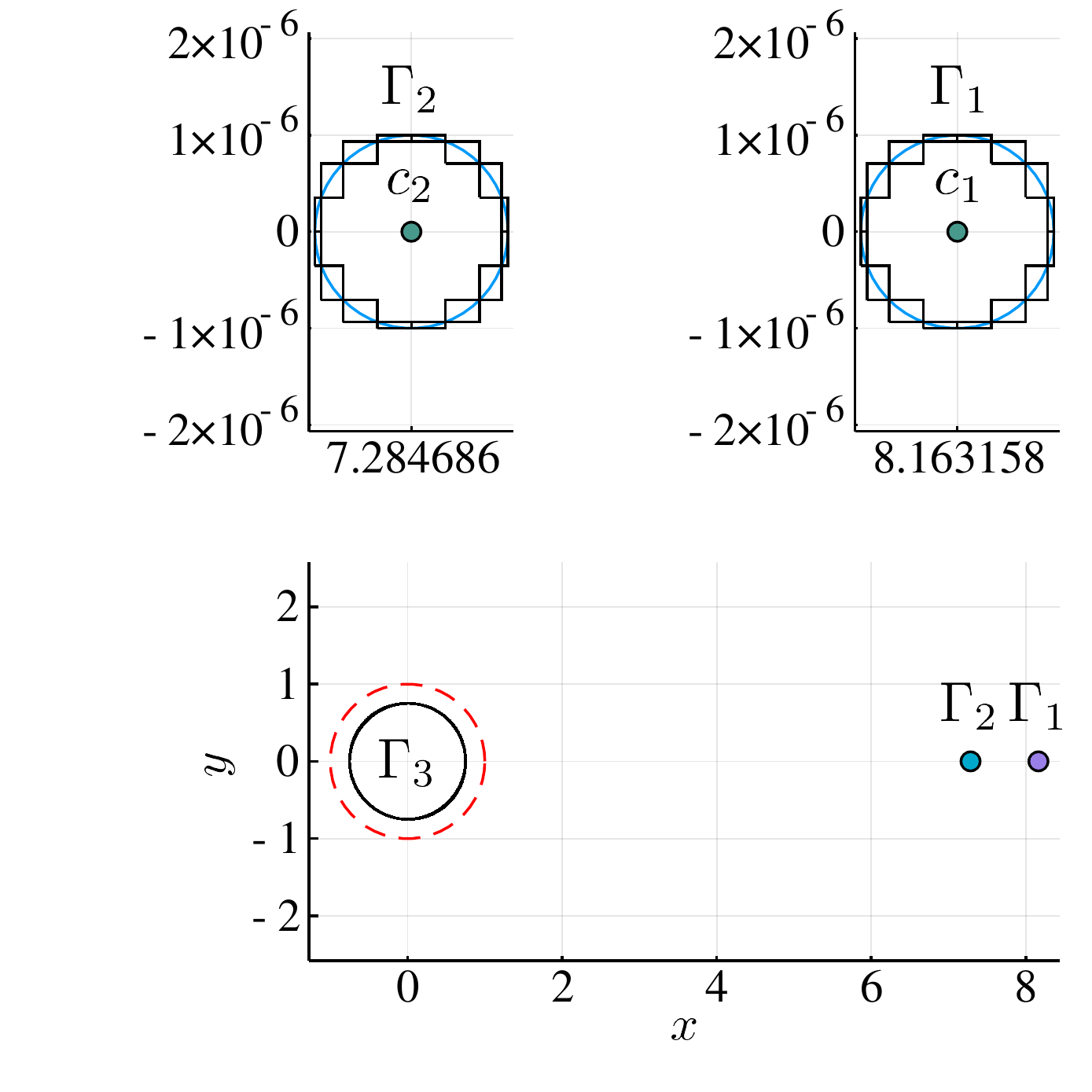}
\caption{Rigorous coverings of the circles
  $\Gamma_1,\Gamma_2,\Gamma_3$ (by $16$, $16$, and $1000$ rectangles,
  respectively) used to bound the determinant
  $\mathrm{det}(L_\mu-\lambda I)$ away from zero for all $\mu\in[0,1]$ and
  $\lambda\in\Gamma_{1,2,3}$ and hence establish that no
  eigenvalues of $L_\mu$ may intersect $\Gamma_{1,2,3}$.  The unit
  circle is shown (dashed) for comparison.\label{fig:spectrum}}
\end{center}\end{figure}

We then consider a smooth one-parameter family of linear operators
\(\mu\mapsto L_\mu\) with \(L_1 = M\) and \(L_0=D\), a diagonal
operator whose spectrum can therefore be determined trivially to have
the correct form. We may then identify disjoint circles
\(\Gamma_1,\Gamma_2,\Gamma_3\) chosen so that \(\Gamma_1,\Gamma_2\)
surround the expanding eigenvalues \(\alpha^4\) and \(\delta\)
respectively, while \(\Gamma_3\) surrounds the rest of the spectrum
within the interior of the unit disc; see Fig. \ref{fig:spectrum}. We
note that the determinant is continuous in the linear operator and, by
proving that \(\mathrm{det}(L_\mu-\lambda I)\) is bounded strictly
away from zero for all \(\mu\in[0,1]\) and all \(\lambda\) on each
circle, we establish that no eigenvalue crosses the circles
\(\Gamma_1,\Gamma_2\), and \(\Gamma_3\). Thus the spectrum of
\(DT(G^{*})\) has the same structure as that of \(D\), with exactly
one eigenvalue bounded within each of \(\Gamma_1\) and \(\Gamma_2\),
and the rest of the spectrum bounded by \(\Gamma_3\)
\cite{kato2013perturbation}.

\subsection{Bounding eigenfunctions and their eigenvalues}\label{bounding-eigenfunctions-and-their-eigenvalues}

Next, we use a novel method to find tight rigorous bounds on
eigenfunction-eigenvalue pairs \((V, \lambda)\) by first rewriting the
eigenproblem in a modified nonlinear form and then adapting the method
used to bound the fixed point \(G^{*}\).  The eigenproblem is given by
\[
(DT(G^{*})-\lambda I)V = 0.
\]
Sticking with the sequence of monomials (expanded with respect to
\(\Omega\)) as Schauder basis, we take \(k\) to be the coordinate
index of the first nonzero coefficient of the eigenfunction
corresponding to \(\delta\) (resp. \(\alpha^4\)), and define
\(\varphi\) to be the corresponding linear coordinate functional. We
choose a normalisation for the eigenfunctions that fixes the
corresponding eigenvalue as the coefficent \(a_k\) of \(V\),
\begin{eqnarray*}
(V, \lambda) &\mapsto& \lambda\frac{V}{\varphi(V)},
\end{eqnarray*}
and rewrite the eigenproblem in a novel way as the corresponding
(nonlinear in \(V\)) problem
\[
F(V) := (DT(G^{*})-\varphi(V))V = 0.
\]

An initial guess, \(V^0\), for the eigenfunction \(V\) may be found by
computing the correponding normalised eigenvector for the truncated
problem nonrigorously and then employing a nonrigorous Newton
iteration to improve the initial guess.

\subsection{Newton's method for eigenfunctions}\label{newtons-method-for-eigenfunctions}

Following the method used to bound the renormalisation fixed point, we
then form a quasi-Newton operator, \(\Psi\), whose fixed points are
the relevant zeros. We first note that \(F\) has Fr\'echet derivative
given formally by
\begin{eqnarray*}
DF(V)\delta V
&=&
DT(G^{*})\delta V - \varphi(\delta V)V - \varphi(V)\delta V.
\end{eqnarray*}
The quasi-Newton operator for this problem is given by
\[
\Psi:V \mapsto V - \widehat\Lambda\left[DT(G^{*})V - \varphi(V)V\right],
\]
in which we choose a fixed invertible linear operator
\(\widehat\Lambda\) such that for all \(f\in B^3:=
B(V^{0};0,\widehat{\rho})\). {We take}
\[
\widehat\Lambda\delta V\simeq\left[DT(G^{*})\delta V-\varphi(\delta V)V^0-\varphi(V^0)\delta V\right]^{-1}.
\]
The Fr\'echet derivative of the quasi-Newton operator is thus given by
\begin{eqnarray*}
D\Psi(V)\delta V 
&=& \delta V - \widehat\Lambda\bigl[D{T}(G^{*})\delta V\nonumber\\
&&\quad{}- \phi(\delta V)V - \phi(V)\delta V\bigr].
\end{eqnarray*}

\subsubsection{Choosing the fixed linear operator}\label{choosing-the-fixed-linear-operator}

Following sections \ref{sec:epsilon} and \ref{sec:kappa}, we aim to
bound \(\|\Psi(V^0)-V^0\|\le\widehat\epsilon\) via function ball
operations on a singleton ball \(B^2:= B(V^0;0,0)\). We must then
bound \(\|D\Psi(V)(e_j)\|\le\kappa<1\) for all \(V\in
B^3=B(V^0;0,\widehat\rho)\) and all \(j\ge 0\).

Anticipating a dependency problem of the sort encountered in section
\ref{sec:depend}, we examine the linear operator, \(\widehat\Lambda\)
more closely. We have
\begin{eqnarray*}
DF(V)\delta V
&=& D{T}(G)\delta V - \varphi(\delta V)V - \varphi(V)\delta V\\
&=& \left(D{T}(G) - Ve_k^{*} - V_kI\right)\delta V,\\
DF(V)
&\simeq& \Delta - V^0e_k^{*} - V^0_kI,
\end{eqnarray*}
where \(V^0\) is a suitable approximate eigenfunction and \(e_k^{*}\)
denotes the adjoint of the basis element \(e_k\), and the subscript on
\(V\) and \(V^0\) denotes the relevant power series
coefficient. Recall that \(\Delta\simeq DT(G^0)\) is chosen so that
its action on \(H\mathscr{A}(\Omega)\) is zero. In order to implement
\(\Lambda\) (which we choose to be the inverse of the above operator)
we need to think about the action of the operator on the polynomial
and high-order parts of the space.

Assume, without loss of generality, that \(k=0\) so that
\(\varphi(V)=V_0\) then, for a suitable \(V^0\) (chosen with
\(HV^0=0\)), we may then take the (block diagonal) operator specified
by
\begin{widetext}
    \[
\Gamma
=\Delta-V^0e_0^{*}-V^0_0I
=
\left(\begin{array}{cccc|c}
\Delta_{00}-2V^0_0 & \Delta_{01} & \cdots & \Delta_{0N} & 0\\
\Delta_{10}-V^0_1 & \Delta_{11}-V^0_0 & \cdots & \Delta_{1N} & 0\\
\vdots & \vdots & \ddots & \vdots & \vdots\\
\Delta_{N0}-V^0_N & \Delta_{N1} & \cdots & \Delta_{NN}-V^0_0 & 0\\
\hline
0 & 0 & \cdots & 0 & -V^0_0I
\end{array}\right).
\]
\end{widetext}

\subsubsection{Overcoming the dependency problem}\label{overcoming-the-dependency-problem}\label{sec:deltadepend}

Recall that
\begin{eqnarray*}
\Psi:V
&\mapsto& V-\widehat\Lambda\bigl[DT(G)V - \varphi(V)V\bigr],
\end{eqnarray*}
with Fr\'echet derivative
\begin{eqnarray*}
D\Psi(V):\delta V
&\mapsto& \delta V - \widehat\Lambda\bigl[DT(G)\delta V\\
&&\quad{}- \varphi(\delta V)V - \varphi(V)\delta V\bigr].
\end{eqnarray*}
We recall that multiple occurences of the perturbation \(\delta V\) in
an expression are treated as functions varying independently within
function balls in the rigorous computational framework, each
contributing separately to the resulting norm. There is therefore a
dependency problem due to the terms \(\delta V\) and
\(\widehat\Lambda\varphi(V)\delta V\) in the above.

To resolve this, consider the action of \(D\Psi(V)\) on a high-order
perturbation \(\delta V_H\in H\mathscr{A}(\Omega)\):
\begin{eqnarray*}
D\Psi(V)\delta V_H
&=&
\delta V_H - \widehat\Lambda\bigl[DT(G)\delta V_H - \varphi(V)\delta V_H\bigr]\nonumber\\
&=&\left(1-\frac{\varphi(V)}{\varphi(V^0)}\right)\delta V_H - \widehat\Lambda DT(G)\delta V_H,
\end{eqnarray*}
since \(\varphi(\delta V_H)=0\) and the action of \(\widehat\Lambda\)
on the high-order part of the space is given by \(-(1/V^0_0)I\). Note
that for \(V\) close to \(V^0\), the contribution from the first term
in the above expression is close to zero.

In order to avoid a bound on \(\|D\Psi(B(V^0;0,\hat\rho))(E_{H})\|\)
exceeding \(2\), we therefore use the latter expression given above
for \(D\Psi(V)\delta V_H\), with \(V\) ranging over the ball
\(B(V^0;0,\widehat\rho)\), when computing \(D\Psi(V)E_{H}\).

Using the parameters with the lowest truncation degree sufficient to
prove the existence of the fixed point \(G^{*}\), given in the first
row of Table \ref{tab:bounds}, we obtain a rigorous bound
\(\|\Psi(V^0)-V^0\|<\widehat{\varepsilon}=3.17\times 10^{-16}\), then
choosing \(\widehat{\rho}=10^{-15}\) yields
\(\|D\Psi(B(V^0;0,\widehat{\rho}))\|<\widehat{\kappa}=1.17\times
10^{-3}\), which establishes that \(\Psi\) is indeed a contraction
mapping on \(B(V^0;0,\widehat{\rho})\). {These crude bounds establish
  that the} eigenvalue satisfies
\(\delta\in[7.28468621706,7.28468621709]\). We use high precision and
high truncation degree to obtain much tighter rigorous bounds on both
the eigenvalue and on the coefficients of the corresponding
eigenfunction \(V^{*}\), as shown in Table \ref{tab:bounds} and Table
\ref{tab:digits}, and are able to prove $325$ significant digits of
$\delta$ correct as a result.

\subsection{Evaluating the eigenfunction on larger domains}\label{evaluating-the-eigenfunction-on-larger-intervals}

We note that the eigenfunction \(V\) satisfies the equation
\begin{widetext}
   
\begin{eqnarray*}
V(X)
&=& \delta^{-1}DT(G)V(X)\\
&=& \delta^{-1}\bigl[-a^{-2}V(1) \cdot G(Q(G(Q(a)X)))\\
&&\quad{}+ a^{-1}\cdot V(Q(G(Q(a)X)))\\
&&\qquad{}+ a^{-1}\cdot G'(Q(G(Q(a)X)))\cdot Q'(G(Q(a)X))\cdot V(Q(a)X)\\
&&\qquad\quad{}+ a^{-1}\cdot G'(Q(G(Q(a)X)))\cdot Q'(G(Q(a)X))\cdot G'(Q(a)X)\cdot Q'(a)V(1)\cdot X\bigr],
\end{eqnarray*}
\end{widetext}
where \(a=G(1)\) and \(\delta=\varphi(V)\). This allows us to evaluate
the eigenfunction \(V(X)\) of \(DT(G^{*})\), and hence
\(v(x)=V(Q(x))\), the corresponding eigenfunction for \(DR(g^{*})\),
over larger subintervals of the real line, by constructing recurrence
relations that utilise the function balls \(B^1\ni G^{*}\) and
\(B^3\ni V\), already computed, as a base case.

\begin{figure}[ht]\begin{center}
\includegraphics[width=0.425\textwidth]{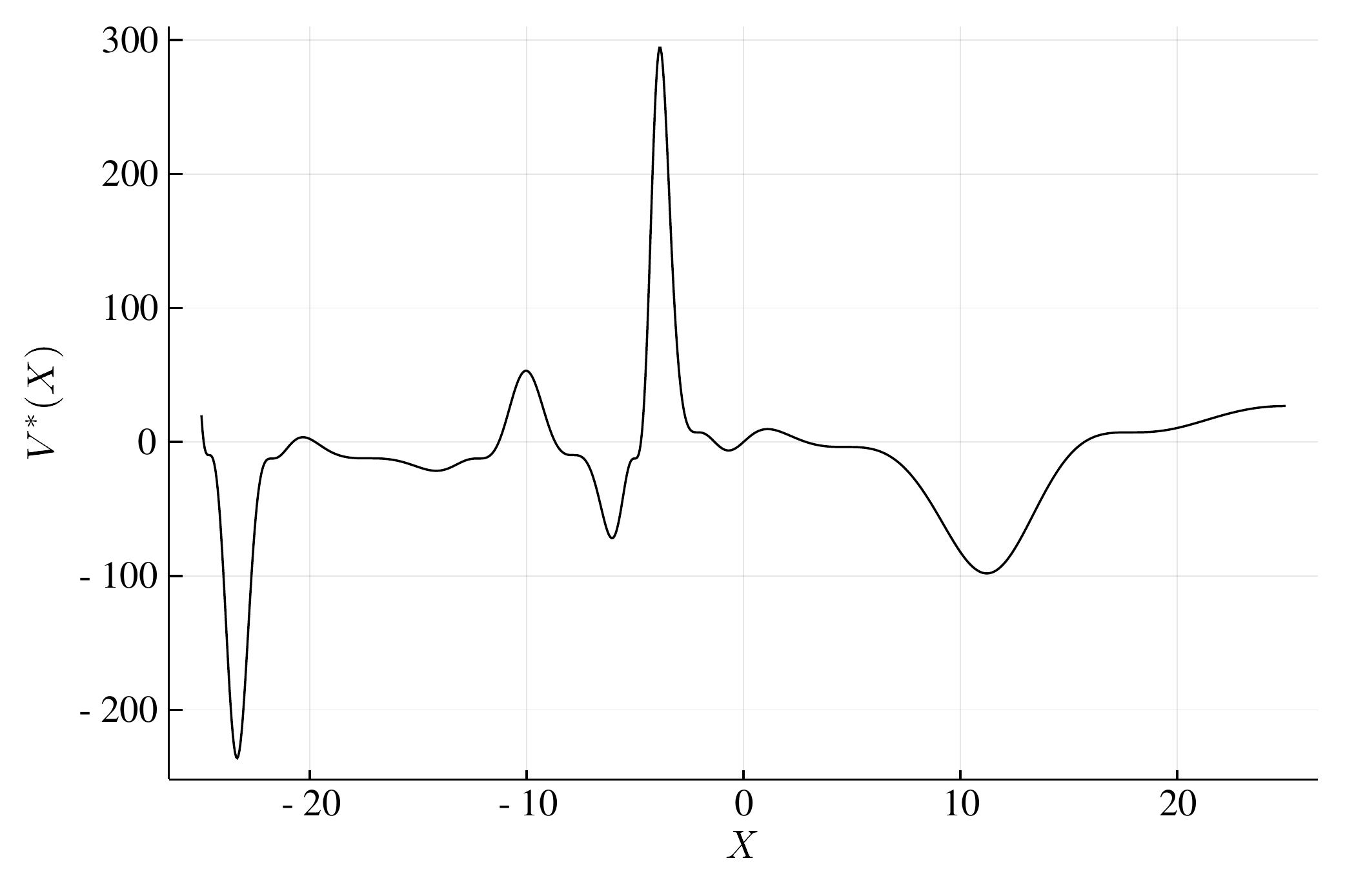}\\
\includegraphics[width=0.425\textwidth]{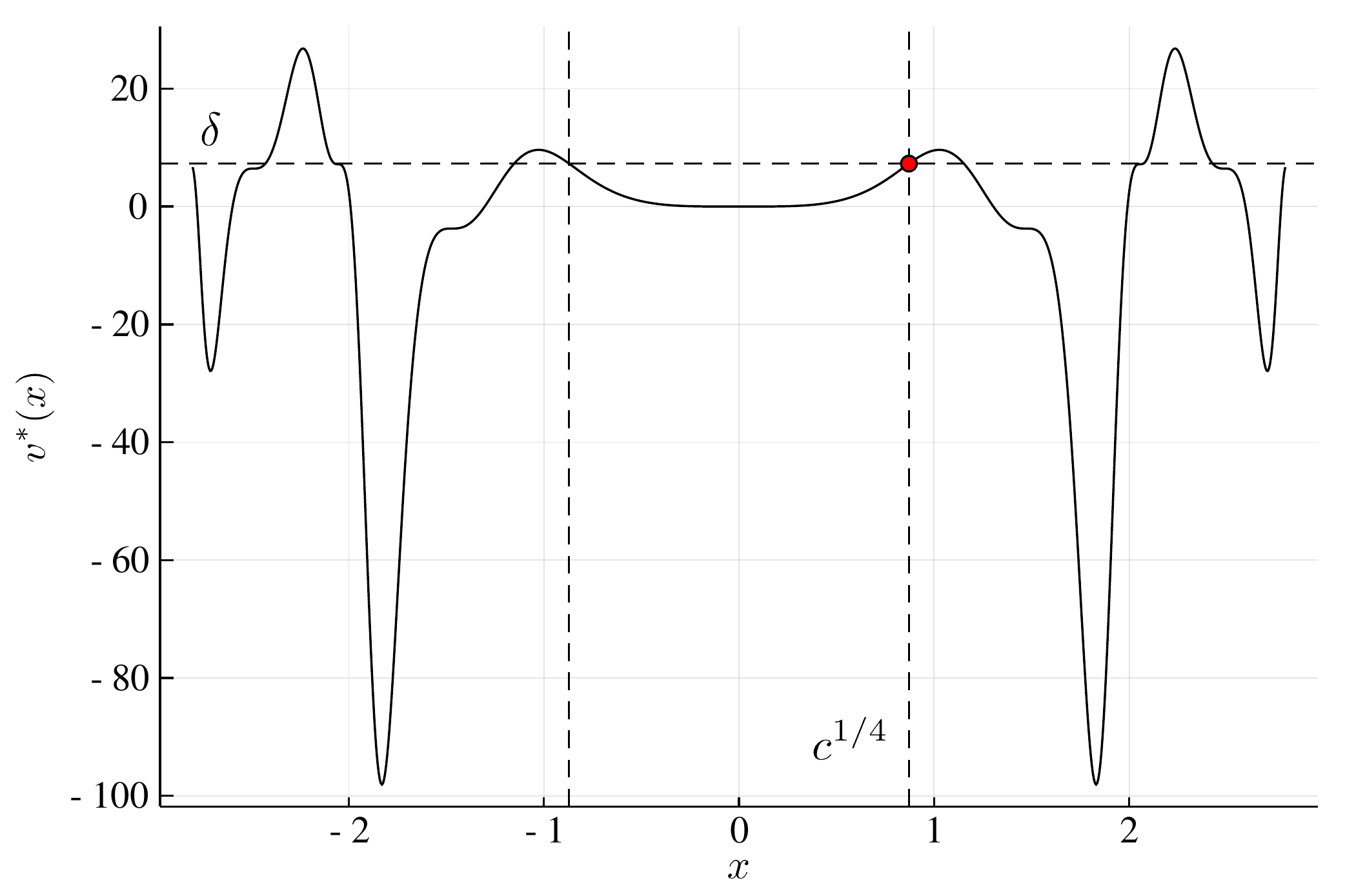}
\caption{(Top) The eigenfunction $V$ corresponding to the essential
  expanding eigenvalue $\delta$.  (Bottom) $v(x) = V(Q(x))$.  The
  dashed lines indicate how the eigenvalue $\delta$ is encoded via the
  chosen normalisation. Here, $\delta=v(c^{1/4})$ since $\delta=V(c)$
  where $\Omega=D(c,r)$; thus $\delta$ corresponds to the power series
  coefficient $a_0$ when $V$ is expanded with respect to $\Omega$.
\label{fig:delta}}
\end{center}\end{figure}

Specifically, we first make use of the fixed-point equation in order
to bound \(G'\) over larger domains: let \(G=G^{*}\), then
\begin{equation}
G(X) = a^{-1}G(Q(G(Q(a)X))).
\end{equation}
Differentiating gives
\begin{eqnarray*}
G'(X)
&=& a^{-1}G'(Q(G(Q(a)X)))\cdot Q'(G(Q(a)X))\\
&&\quad{}\cdot G'(Q(a)X)\cdot Q(a).
\end{eqnarray*}
Using the above expression (together with the fixed-point equation for
\(G\)) recursively allows us to bound \(G'\) and hence, in combination
with the above, \(V\), over larger intervals extending outside
\(\Omega\cap\mathbb{R}\).  See Fig. \ref{fig:delta}.

\section{Critical scaling of added noise}\label{critical-scaling-of-additive-noise}\label{sec:noise}

We now find tight rigorous bounds on the eigenfunction, \(w\), and
eigenvalue, \(\gamma\), controlling the universal scaling of additive
uncorrelated noise. The iteration of a prototypical one-parameter
family, \(f_\mu\), is modified to give \(x_{n+1}=F_{\mu, n}(x_n):=
f_\mu(x_n)+\varepsilon\xi_n\) where, in the simplest case, the
\(\xi_n\) are i.i.d. random variables, independent of the
\(x_n\). Adapting the arguments presented in
\cite{crutchfield1981scaling,shraiman1981scaling}, and retaining the
deterministic scaling \(a=f(1)\) in the definition of the
renormalisation operator, we write \(w=W\circ Q\) and consider the
{modified} eigenproblem
\[
\gamma^2W = \mathcal{L}W,
\]
in which we define the linear operator \(\mathcal{L}\) by
\[
\mathcal{L}W := {L_1}^2\cdot W(Q(G(Q(a)X))) + {L_2}^2\cdot W(Q(a)X),
\]
where we define
\begin{eqnarray*}
L_1 &:=& a^{-1},\\
L_2 &:=& a^{-1}G'(Q(G(Q(a)X)))\cdot Q'(G(Q(a)X)).
\end{eqnarray*}
In the above, \(G=G^{*}\). We note that the expressions \(L_1,L_2\)
are those prefactors {of $\delta G(\cdot)$} in the Fr\'echet
derivative, \(DT(G)\), of equation \ref{eqn:frechet} that do not
correspond to variations in \(a\) (equivalently, those in the terms
\ref{eqn:drg2},\ref{eqn:drg3} of \(DR(g)\)). We note also that the
corresponding operator acting on \(g\) emerges as a special case of
the analysis presented in \cite{fiel1987correlated} for the correlated
case. Following our treatment for the eigenfunctions of \(DT(G)\), we
{take the novel approach of encoding} the eigenvalue within \(W\) by
defining \(\gamma=\varphi(W)\) and expressing the problem in the
modified nonlinear form
\[
\mathcal{F}(W)
:= \left(\mathcal{L}-\varphi(W)^2I\right)W
= 0.
\]

\begin{figure}[ht]\begin{center}
\includegraphics[width=0.425\textwidth]{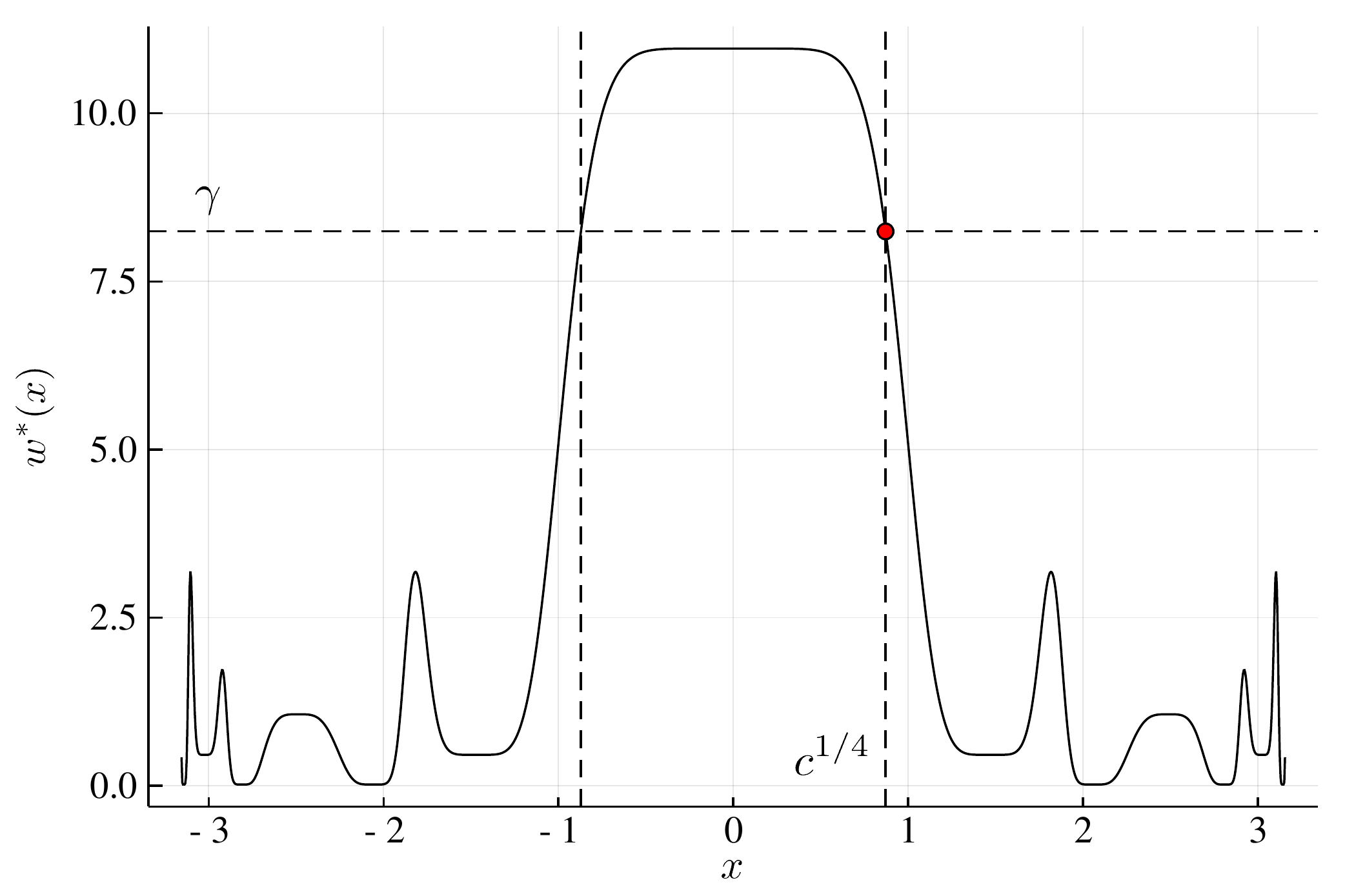}\\
\includegraphics[width=0.425\textwidth]{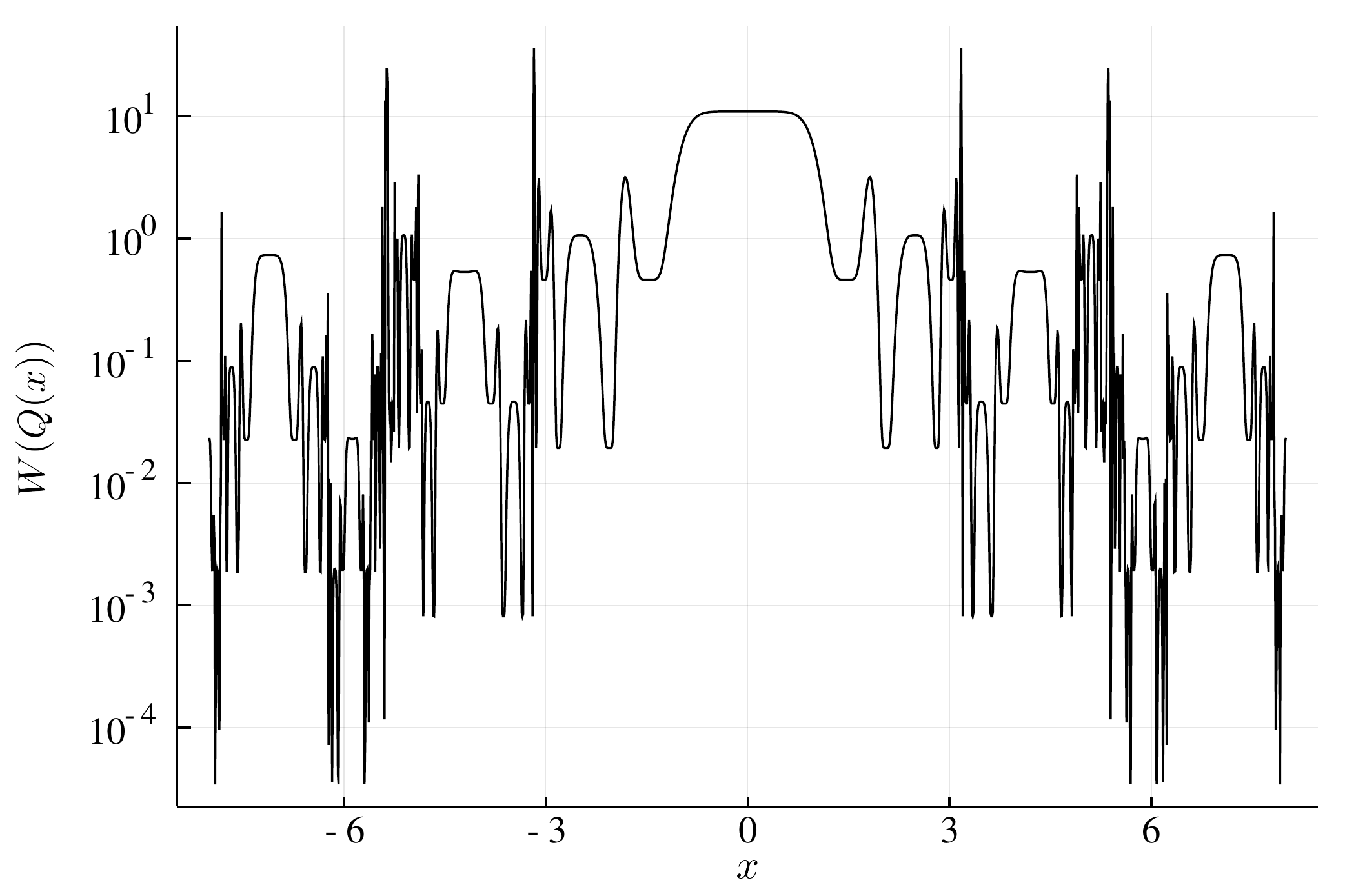}
\caption{(Top) The eigenfunction $w(x)=W(Q(x))$ corresponding to
  critical scaling of Gaussian noise in the iteration of maps with
  quartic critical point. The dotted lines indicate how $\gamma$ is
  encoded via the chosen normalisation. Here,
  $\gamma=W(c)=w(c^{1/4})$.  (Bottom) Plotting $w(x)$ with a
  logarithmic scale over a larger interval emphasises the self-similar
  structure.  (Observe that $w(x)\ge 0$ for all $x$.)\label{fig:noise}}
\end{center}\end{figure}

The operator \(\mathcal{F}\) has Fr\'echet derivative
\begin{eqnarray*}
D\mathcal{F}(W):\delta W
&\mapsto& \mathcal{L}\delta W - 2\varphi(W)\varphi(\delta W)W - \varphi(W)^2\delta W.
\end{eqnarray*}

We form the quasi-Newton operator
\[
\Theta(W) := W-\Lambda\mathcal{F}(W),
\]
where \(\Lambda\) is a fixed linear operator
\(\Lambda\simeq[D\mathcal{F}(W^0)]^{-1}\). The Fr\'echet derivative is
given by
\begin{eqnarray*}
&&D\Theta(W)\delta W\\
&=& \delta W-\Lambda D\mathcal{F}(W)\delta W\\
&=& \delta W-\Lambda\bigl[\mathcal{L}\delta W - 2\varphi(W)\varphi(\delta W)W - \varphi(W)^2\delta W\bigr].
\end{eqnarray*}
In particular, we take
\[
D\mathcal{F}(W)\simeq\mathcal{L}-2\varphi(W^0)W^0e_0^{*}-\varphi(W^0)^2I,
\]
choosing \(W^0\) such that \(HW^0=0\), and take \(\Lambda\) to be the
inverse operator corresponding to the right-hand side, which therefore
has the following action on high-order terms, {$\delta W_H\in HA$,}
\[
\Lambda\delta W_H = -\frac{1}{\varphi(W^0)^2}\delta W_H.
\]

To mitigate the corresponding dependency problem, we compute the
action of \(D\Theta(W)\) on a high-order perturbation \(\delta W_H\in
H\mathscr{A}(\Omega)\):
\[
D\Theta(W)\delta W_H
=\left[1-\left(\frac{\varphi(W)}{\varphi(W^0)}\right)^2\right]\delta W_H
- \Lambda\mathcal{L}\delta W_H,
\]
noting again that for $W$ close to $W^0$ the first term is close to
zero.

Using the parameters for \(G^{*}\) given in the first row of Table
\ref{tab:bounds}, we obtain
\(\|\Theta(W^0)-W^0\|<\tilde{\varepsilon}=2.35\times 10^{-16}\);
choosing \(\tilde{\rho}=10^{-15}\) then gives
\(\|D\Theta(B(W^0;0,\tilde{\rho}))\|<\tilde{\kappa}=7.85\times
10^{-3}\), establishing that \(\Theta\) is a contraction on
\(B(W^0;0,\tilde{\rho})\). Table \ref{tab:bounds} demonstrates that
these bounds may be improved significantly. Fig. \ref{fig:noise} shows
the corresponding eigenfunction. Working with truncation degree \(40\)
and \(40\) digits in the significand yields the crude bound
\(\gamma\in[8.24391085424, 8.24391085427]\) for the noise eigenvalue,
which we again improve by taking higher truncation degree and by using
multiprecision arithmetic. Table \ref{tab:digits} shows 323 digits
proven correct. These bounds confirm the initial digits of the
numerical estimate presented in \cite{kuznetsov2002generalized}.

\section{Conclusions}\label{conclusions}

We have obtained tight bounds on the renormalisation fixed point
function for period doubling in unimodal maps with critical point of
degree \(4\), by means of a rigorous computer-assisted existence proof
using the contraction mapping theorem on a suitable space of analytic
functions. We have established the structure of the spectrum of the
linearised operator at the fixed point, providing bounds on expanding
eigenvalues. By expressing the corresponding eigenproblem in nonlinear
form, we have used a contraction mapping argument to provide rigorous
bounds on eigenfunction-eigenvalue pairs, and have adapted the
technique to bound the eigenfunction and eigenvalue controlling the
universal scaling of added noise in the case of a deterministic choice
of normalisation in the renormalisation operator. These techniques
deliver tight bounds on the relevant analytic functions and the
corresponding universal constants.

The method may be adapted to unimodal maps with general integer
critical exponent. In the case of general even degree critical points,
this relies on finding suitable function domains. In the case of odd
degree critical points, the method may also be applied by recourse to
a suitably-modified functional equation. Increasing the degree will
inevitably lead to challenges in the rigorous numerics. We examine
both cases in forthcoming publications and use the bounds thus
obtained to gain rigorous bounds on the Hausdorff dimensions of the
relevant attractors at the accumulation of the period-doubling
cascade.

% We use these bounds to gain bounds on the Hausdorff dimension of
% the corresponding attractor at the limit of the period-doubling
% cascade for families in the relevant universality class.

\subsection{Computational issues}\label{sec:computational-issues}

In the above computations, the tightness of the bound on the
contractivity, \(\kappa\), of the three quasi-Newton methods is
limited by the high-order bound on \(\|D\Phi(B)E_{H}\|\) for the
relevant operator $\Phi$ and function ball $B$.  Recall that the ball
$E_H$ is the convex hull of all high-order basis elements.  The bound
computed on the supremum of the quantities $\|D\Phi(G)e_k\|$ is
dominated by this high-order contribution for \(N=40,80,160,320\).
However, for \(N=480,640\), the bound on the supremum is achieved by
one of the \(D\Phi(B^{1})E_k\) for \(0\le k\le N\), indicating that a
sufficently high truncation degree has been taken such that the loss
of information concerning the distribution of the high-order bound
amongst high-order coefficients no longer provides the dominant
obstacle to improving the bound on contractivity.

The computations were verified independently by two different
implementations of the function ball algebra: the first is written in
the high-performance language Julia \cite{bezanson2017julia} and
utilises multi-precision binary floating-point arithmetic with
rigorous directed rounding modes \cite{kaucher2014self}.  The second
is written in the language Python and utilises multi-precision decimal
floating-point arithmetic with rigorous directed rounding modes.  The
binary arithmetic conforms to the relevant subset of standard
IEEE754-2008 and the decimal arithmetic to the relevants subsets of
standards ANSI X3.274-1996, IEEE754-2008, and ISO/IEC/IEEE60559:2011.

The framework for rigorous function ball operations is adapted from
that of \cite{ew1985computer,eckmann1982existence}, specialised to the
disc algebra $A$, and implemented with multi-precision
arithmetic. Optimisations were made for the computation of the
high-order bound on products. Parallel computation was used for
obtaining the bounds on contractivity over the balls containing the
basis elements $e_k$. Closures were used carefully in order to avoid
recomputation of bounds on common sub-expressions in the corresponding
Fr\'echet derivatives.

The integrity of the frameworks is verified with the aid of over
\(1200\) unit tests and functional tests. Where parallel computation
has been used, care was taken to use multiprocessing rather than
threads in order to protect the integrity of rounding modes across
processes and the results of parallel computations were verified
against the corresponding serial code.
\begin{acknowledgments}
We thank Andreas Stirnemann and Ben Mestel for helpful discussions. 
\end{acknowledgments}\vfil
    
\bibliography{references}
    
\end{document}